\documentclass[10pt]{amsart}

\usepackage{a4}
\usepackage{amssymb}

\newtheorem{thm}{Theorem}[section]
\newtheorem{lem}[thm]{Lemma}

\theoremstyle{definition}
\newtheorem{dfn}[thm]{Definition}
\newtheorem{ntn}[thm]{Notation} 
\newtheorem{rmk}[thm]{Remark}
\newtheorem{exa}[thm]{Example}
\newtheorem{con}[thm]{Convention}
\newtheorem{ass}[thm]{Assumption}

\newenvironment{pf}
  {\noindent{{\em Proof: }}}
  {\qed \medskip}

\def\quod{\hskip 0.2em\relax }

\newcommand{\abs}[1]{\lvert#1\rvert}
\newcommand{\Pb}{\mathbf{P}}

\newcommand{\Q}{\mathbb{Q}}
\newcommand{\Z}{\mathbb{Z}}

\newcommand{\Aut}{\mathrm{Aut}}
\newcommand{\s}{\mathbb{S}}
\newcommand{\Ql}{\mathbb{Q}_{\ell}}
\newcommand{\Vla}{\mathbb{V}_{\lambda}}

\newcommand{\F}{\mathbb{F}}

\newcommand{\M}[1]{\mathcal M_{{#1}}}
\newcommand{\Mm}[2]{\mathcal M_{{#1},{#2}}}
\newcommand{\Hh}[1]{\mathcal H_{{#1}}}
\newcommand{\Hhm}[2]{\mathcal H_{{#1},{#2}}}

\newcommand{\Mmb}[2]{\overline{\mathcal M}_{{#1},{#2}}}

\newcommand{\Ll}{\mathbf L}

\newcommand{\Eulc}{\mathbf{e}_c}

\newcommand{\PGLk}[1]{\abs{\mathrm{PGL}_{#1}(k)}}
\newcommand{\linsub}[1]{U_{#1,\lambda}}
\newcommand{\zee}[1]{z(#1)}

\newcommand{\Pbk}{\mathbf{P}_k}
\newcommand{\Pbpl}{\abs{\mathbf{P}^2_k(\lambda)}}

\newcommand{\Qq}[1]{\mathcal{Q}_{{#1}}}
\newcommand{\Qqt}[1]{\tilde{\mathcal{Q}}_{{#1}}}
\newcommand{\Mmt}[2]{\tilde{\mathcal M}_{{#1},{#2}}}

\begin{document}
\pagestyle{plain}

\title[Cohomology of moduli spaces of curves of genus $3$ via point counts]{Cohomology of moduli spaces of curves of genus three via point counts}  
\author{Jonas Bergstr\"om}
\email{jonasb@math.kth.se}
\address{Department of Mathematics, KTH, S--100 44 Stockholm, Sweden}

\begin{abstract}
In this article we consider the moduli space of smooth $n$-pointed non-hyperelliptic curves of genus $3$. In the pursuit of cohomological information about this space, we make $\mathbb{S}_n$-equivariant counts of its numbers of points defined over finite fields for $n \leq 7$. 
Combining this with results on the moduli spaces of smooth pointed curves of genus $0$, $1$ and $2$, and the moduli space of smooth hyperelliptic curves of genus $3$, we can determine the $\mathbb{S}_n$-equivariant Galois and Hodge structure of the ($\ell$-adic respectively Betti) cohomology of the moduli space of stable curves of genus $3$ for $n \leq 5$ (to obtain $n \leq 7$ we would need counts of ``$8$-pointed curves of genus $2$''). 
\end{abstract}

\maketitle

\section{Introduction} \label{sec-intr}
Inside the moduli space $\Mm gn$ of $n$-pointed smooth curves of genus $g$ there is a closed subset $\Hhm gn$ consisting of the hyperelliptic curves. In the case of genus $3$ we denote the complement of this closed subset by $\Qq{n}$. In this article we will compute the number of points over finite fields of this space for all $n \leq 7$. On $\Mm gn$ we have an action of the symmetric group $\s_n$ by permuting the $n$ marked points on the curves. To take this action into account we will be counting $\s_n$-equivariantly.

There is cohomological information to be found from counting points over finite fields. For instance, if the number of points over any finite field $\F_q$ of a smooth and proper DM-stack is a polynomial in $q$, then the Galois and Hodge structures of the ($\ell$-adic respectively Betti) cohomology groups of the stack are determined by this polynomial, see Theorem 2.1 in \cite{EB} by van den Bogaart-Edixhoven. 

The moduli space $\Mmb gn$ of $n$-pointed stable curves of genus $g$ is a proper and smooth DM-stack which contains $\Mm gn$ as an open part. For all $n \leq 5$, the $\s_n$-equivariant count of points of $\Qq{n}$ is found in Section \ref{sec-stable} to be the only missing result to make an $\s_n$-equivariant count of points of $\Mmb{3}{n}$. These counts are all found to be polynomial (this is also true for $\Qq{6}$ and $\Qq{7}$) and hence we can apply a version of the theorem above to $\Mmb{3}{n}$ for $n \leq 5$. This theorem also gives the $\s_n$-equivariant Galois and Hodge Euler characteristic of $\Mm{3}{n}$ for $n \leq 5$. In Section \ref{sec-localsystems} we present our results for $n=6$ and $7$, formulated in terms of Galois Euler characteristics of some natural local systems on $\M{3}$. The results for $n=0$ and $1$ are the only previously known ones, by the work of Looijenga in \cite{Looij}, and Getzler-Looijenga in \cite{GL}. In the thesis \cite{Orsolathesis} by Tommasi, the result for $\Mm{3}{2}$ is proved by a different method. 

The elements of $\Qq{n}$ can be embedded into the projective plane as non-singular quartic curves using their canonical linear system. An $\s_n$-equivariant count of points over a finite field $k$ of $\Qq n$ can then be achieved by finding, for each $n$-tuple of points $(p_1,\ldots,p_n)$ in the plane that as a set is defined over $k$, the number of non-singular plane quartics over $k$ that pass through these points, see Section \ref{sec-quartics}. 
The method we will use to compute these numbers is presented in Section \ref{sec-sieve}. We begin with the linear subspace, of the $\Pb^{14}(k)$ of all plane quartic curves defined over $k$, consisting of the curves that pass through the points $(p_1,\ldots,p_n)$. We then add and subtract subspaces of singular quartic curves that contain $(p_1,\ldots,p_n)$ using the sieve method. That is, in such a way that in each step, more and more singular curves will have been removed precisely once. This procedure must be stopped at some point, due to the quartics with infinitely many singularities. We will choose to do this already after one step, see Section \ref{sec-M1}. All quartics with precisely one singularity will then have been removed once. In Section \ref{sec-linsub} we compute the dimensions of the linear subspaces that appear in this procedure. We then have to amend for the curves with two or more singularities so that they also will have been removed precisely once. The information necessary to do this is found in Sections \ref{sec-library} to \ref{sec-two}, together with Section~\ref{sec-cubics}. In the latter section we apply the same method as for the quartics to the cubics, giving $\s_n$-equivariant counts of points of $\Mm{1}{n}$ for $n \leq 10$. These results are all in agreement with the work of Getzler \cite{G-res}, where the $\s_n$-equivariant Hodge Euler characteristic of $\Mm{1}{n}$ is determined for any $n$.

\section*{Acknowledgements}
I thank Torsten Ekedahl for help with Theorem \ref{thm-trgal} and Institut Mittag-Leffler for support during the preparation of this paper. I am grateful to my advisor Carel Faber for all help. 

\section{Cohomology of the moduli space of stable curves of genus $3$} \label{sec-stable}
There is a stratification of $\Mmb gn$ into pieces that are related to $\Mm{\tilde g}{\tilde n}$ for $\tilde g \leq g$ and $\tilde n \leq n+2(g-\tilde g)$. Beginning with $\s_{\tilde n}$-equivariant counts of points defined over a finite field (see Section \ref{sec-quartics}) of all these moduli spaces $\Mm{\tilde g}{\tilde n}$, we can get an $\s_{n}$-equivariant count of points of $\Mmb gn$ using this stratification. If these counts, when considered as functions of the number of elements of the finite field, are polynomials, we obtain the $\s_n$-equivariant $\mathrm{Gal}(\bar{\Q}/\Q)$ and Hodge structure of the ($\ell$-adic resp. Betti) cohomology of $\Mmb gn$ by Theorem~3.4 in \cite{Mbar4}. Hence, using the following $\s_n$-equivariant counts of points: 
\begin{itemize}
\item[$\star$] $\Mm{0}{n}$ for $n \leq 11$, 
which is found in \cite{Lehrer};
\item[$\star$] $\Mm{1}{n}$ for $n \leq 9$, 
which is found in Section \ref{sec-cubics}; 
\item[$\star$] $\Mm{2}{n}$ for $n \leq 7$, which is found in \cite{Jonas1};
\item[$\star$] $\Hhm{3}{n}$ for $n \leq 5$, which is found in \cite{Jonas1};
\item[$\star$] $\Qq{n}$ for $n \leq 5$, which is found in this article,
\end{itemize}
which are all polynomial in the sense above, we get these results for $\Mmb {3}{n}$ for $n \leq 5$. By Theorem~3.4 in \cite{Mbar4} we also get the corresponding results for $\Mm {3}{n}$ for $n \leq 5$, but in terms of Euler characteristics.

The piece missing to go one step further is the $\s_8$-equivariant count of points of $\Mm{2}{8}$. The $\s_n$-equivariant computer counts of points over small finite fields, of $\Mm{2}{n}$ for all $n$, of Faber-van der Geer in \cite{Faber-Geer1} and \cite{Faber-Geer2}, show that polynomiality fails for $n=10$, and their conjectural interpretation suggests that it holds for $8$ and $9$. 

In the theorems below we give $\s_n$-equivariant Hodge Euler characteristics (which by purity are sufficient to conclude the Hodge structure) in terms of the Schur polynomials and $\Ll$, the class of the Tate Hodge structure of weight $2$ in the Grothen\-dieck group of rational Hodge structures. The results for $n=0$ and $1$ can be found in \cite{GL}.

\begin{thm}
The equivariant Hodge Euler characteristic of $\Mmb{3}{2}$ is equal to 
\begin{gather*}
(\Ll^8+7\Ll^7+31\Ll^6+74\Ll^5+100\Ll^4+74\Ll^3+31\Ll^2+7\Ll+\mathbf{1})s_{2} \\
+\, \Ll(2\Ll^6+11\Ll^5+30\Ll^4+42\Ll^3+30\Ll^2+11\Ll+\mathbf{2}) s_{1^2}
\end{gather*}
\end{thm}

\begin{thm}
The equivariant Hodge Euler characteristic of $\Mmb{3}{3}$ is equal to 
\begin{gather*}
(\Ll^9+9\Ll^8+50\Ll^7+157\Ll^6+277\Ll^5+277\Ll^4+157\Ll^3+50\Ll^2+9\Ll+\mathbf{1})s_{3}\\
+\, \Ll(4\Ll^7+33\Ll^6+120\Ll^5+228\Ll^4+228\Ll^3+120\Ll^2+33\Ll+\mathbf{4}) s_{21}\\
+\, \Ll^2(2\Ll^5+13\Ll^4+27\Ll^3+27\Ll^2+13\Ll+\mathbf{2})s_{1^3}
\end{gather*}
\end{thm}

\begin{thm}
The equivariant Hodge Euler characteristic of $\Mmb{3}{4}$ is equal to 
\begin{gather*}
(\Ll^{10}+11\Ll^9+76\Ll^8+296\Ll^7+676\Ll^6+887\Ll^5+676\Ll^4+296\Ll^3+76\Ll^2+\ldots )s_{4} \\
+\,\Ll(6\Ll^8+63\Ll^7+305\Ll^6+762\Ll^5+1032\Ll^4+762\Ll^3+305\Ll^2+63\Ll+\mathbf{6}) s_{31}\\
+\,\Ll(2\Ll^8+27\Ll^7+136\Ll^6+351\Ll^5+478\Ll^4+351\Ll^3+136\Ll^2+27\Ll+\mathbf{2}) s_{2^2}\\
+\,\Ll^2(9\Ll^6+69\Ll^5+202\Ll^4+288\Ll^3+202\Ll^2+69\Ll+\mathbf{9})s_{21^2}\\
+\,2\Ll^3(\Ll^4+4\Ll^3+6\Ll^2+4\Ll+\mathbf{1})s_{1^4}
\end{gather*}
\end{thm}

\begin{thm}
The equivariant Hodge Euler characteristic of $\Mmb{3}{5}$ is equal to 
\begin{gather*} 
(\Ll^{11}+13\Ll^{10}+105\Ll^9+500\Ll^8+1419\Ll^7+2397\Ll^6+2397\Ll^5+1419\Ll^4+\ldots )s_{5}\\
+\,\Ll(8\Ll^9+106\Ll^8+637\Ll^7+2043\Ll^6+3633\Ll^5+3633\Ll^4+2043\Ll^3+\ldots )s_{41}\\
+\,\Ll(4\Ll^9+66\Ll^8+444\Ll^7+1501\Ll^6+2734\Ll^5+2734\Ll^4+1501\Ll^3+\ldots )s_{32}\\
+\,\Ll^2(20\Ll^7+202\Ll^6+793\Ll^5+1528\Ll^4+1528\Ll^3+793\Ll^2+202\Ll+\mathbf{20})s_{31^2}\\
+\,2\Ll^2(6\Ll^7+61\Ll^6+244\Ll^5+474\Ll^4+474\Ll^3+244\Ll^2+61\Ll+\mathbf{6})s_{2^21}\\
+\,\Ll^3(13\Ll^5+77\Ll^4+172\Ll^3+172\Ll^2+77\Ll+\mathbf{13})s_{21^3}\\
+\,\Ll^4(\Ll^3+2\Ll^2+2\Ll+\mathbf{1})s_{1^5}
\end{gather*}
\end{thm}

\section{Cohomology of local systems on $\M 3$ and $\Qq{}$} \label{sec-localsystems} 
In this section the results for $n=6$ and $7$ will be given in terms of some natural local systems. Define the local system $\mathbb{V}:=R^1\pi_{*}(\Ql)$, where $\pi : \M{3,1} \to \M{3}$ is the universal curve. For every $\lambda=(\lambda_1\geq\lambda_2\geq\lambda_3 \geq 0)$ we get an induced local system $\mathbb{V}_{\lambda}$ from the irreducible representation of $\mathrm{GSp}(6)$ with highest weight $(\lambda_1-\lambda_2)\gamma_1+(\lambda_2-\lambda_3)\gamma_2+\lambda_3\gamma_3-\abs{\lambda}\eta$, where $\gamma_i$ are suitable fundamental roots and $\eta$ the multiplier representation. Making $\s_{\tilde n}$-equivariant counts of points of $\Mm{3}{\tilde n}$ over a finite field $k$, for all $\tilde n \leq n$, is equivalent to computing the trace of Frobenius on the compactly supported $\ell$-adic Euler characteristic $\Eulc(\M{3} \otimes \bar k,\mathbb{V}_{\lambda})$, for every $\lambda$ with $\abs{\lambda} \leq n$, and where $\ell\nmid \abs{k}$ (for more details see \cite{G-res}). This works equally well for $\Vla'$, the restriction of $\Vla$ to $\Qq{}$, and in the following theorem we apply this to weight $6$ local systems, where we find a dependence upon the characteristic. 

\begin{thm} Let $k$ be a finite field with $\abs{k}=q$, $\Qq{}':=\Qq{} \otimes \bar k$ and $F$ the geometric Frobenius. If we define $\epsilon$ to be $-1$ if $\mathrm{char}(k) = 2$ and $0$ if $\mathrm{char}(k) \neq 2$, then:
\centerline{
\vbox{
\bigskip
\offinterlineskip
\hrule
\halign{&\vrule#& \quod \hfil#\hfil \!\strut  \cr
height2pt&\omit&&\omit&&\omit&&\omit&&\omit&&\omit& \cr 
& $\lambda$ && $\mathrm{Tr}\bigl(F,\Eulc(\Qq{}',\Vla')\bigr)$ && $\lambda$ && $\mathrm{Tr}\bigl(F,\Eulc(\Qq{}',\Vla')\bigr)$ && $\lambda$ && $\mathrm{Tr}\bigl(F,\Eulc(\Qq{}',\Vla')\bigr)$ & \cr
height2pt&\omit&&\omit&&\omit&&\omit&&\omit&&\omit& \cr 
\noalign{\hrule}
height2pt&\omit&&\omit&&\omit&&\omit&&\omit&&\omit& \cr 
&$(6,0,0)$ && $-q^3+q+4+\epsilon$ && $(4,1,1)$ && $q^4-q^3-q^2$ && $(2,2,2)$ && $q^2+q+2+\epsilon$  &\cr
height2pt&\omit&&\omit&&\omit&&\omit&&\omit&&\omit& \cr 
&$(5,1,0)$ && $-2q^2+2q+1$ && $(3,3,0)$  && $q^5+q$ && && &\cr
height2pt&\omit&&\omit&&\omit&&\omit&&\omit&&\omit& \cr 
& $(4,2,0)$ && $-2q^3+2q+3+\epsilon$ && $(3,2,1)$ && $-q^3+q$ && && &\cr
height2pt&\omit&&\omit&&\omit&&\omit&&\omit&&\omit& \cr 
} \hrule}}
\end{thm}

Denote by $K_0(\mathsf{Gal}_{\Q})$ the Grothendieck group of $\mathrm{Gal(\bar{\Q}/\Q})$-modules.

\begin{thm} \label{thm-trgal} Let $X \xrightarrow{\rho} Spec(\Z)$ be a separated scheme of finite type and $\mathcal F$ an $\ell$-adic constructible sheaf on $X$. The trace of Frobenius on $\Eulc(X \otimes \bar{\F}_p,\mathcal F)$, for almost all (i.e., for all but a finite number of) primes $p$, determines $\Eulc(X \otimes \bar{\Q},\mathcal F)$ as an element in $K_0(\mathsf{Gal}_{\Q})$. 
\end{thm} 
\begin{pf} The sheaf $\rho_{!}\mathcal F$ is constructible and commutes with base change (cf. \cite[Thm. VI.3.2]{Milne}), and thus the specialization morphism $H_c^i(X \otimes \bar{\F}_p,\mathcal F) \to H_c^i(X \otimes \bar{\Q}_p,\mathcal F)$ is an isomorphism for almost all primes $p$. Hence, for these primes, $H_c^i(X \otimes \bar{\Q},\mathcal F)$ is an unramified ($\ell$-adic) $\mathrm{Gal(\bar{\Q}/\Q})$-module and thus Frobenius elements are defined up to conjugacy. By continuity and the Chebotarev density theorem (cf. \cite[Prop. 2.6]{fermat}) we know that the trace of Frobenius for almost all primes determines the character of a $\mathrm{Gal(\bar{\Q}/\Q})$-module. The character, in turn, determines the semi-simplification of the module, i.e., as an element in $K_0(\mathsf{Gal}_{\Q})$. Now apply the above to $\Eulc(X \otimes \bar{\Q},\mathcal F)$. 
\end{pf}

Let us apply Theorem \ref{thm-trgal} to the push-down of $\Vla$ to the coarse moduli space of $\M 3$ for weight $6$ and $7$, where we add to $\Qq{}$ the contribution of $\Hh 3$ in weight $6$ (in weight $7$ it is $0$) found in \cite{Jonas1}. Note that the dependence upon characteristic does not survive on $\M{3}$, which is indeed proved in \cite[Cor. 3.2]{jbgvdg}. 

\begin{thm} If $\mathbf{q}$ denotes the class in $K_0(\mathsf{Gal}_{\Q})$ of the cyclotomic character $\Ql(-1)$ then: 

\centerline{
\vbox{
\bigskip
\offinterlineskip
\hrule
\halign{&\vrule#& \quod \hfil#\hfil \!\strut  \cr
height2pt&\omit&&\omit&&\omit&&\omit& \cr 
& $\lambda$ && $\Eulc(\M{3}\otimes \bar{\Q},\Vla)$ && $\lambda$ && $\Eulc(\M{3}\otimes \bar{\Q},\Vla)$  & \cr
height2pt&\omit&&\omit&&\omit&&\omit& \cr 
\noalign{\hrule}
height2pt&\omit&&\omit&&\omit&&\omit& \cr 
&$(6,0,0)$ && $-\mathbf{q}^3$ && 
$(7,0,0)$ && $-\mathbf{q}^3-\mathbf{q}^2$  &\cr
height2pt&\omit&&\omit&&\omit&&\omit& \cr 
&$(5,1,0)$ && $\mathbf{q}^2-\mathbf{2}$ && 
$(6,1,0)$ && $-3\mathbf{q}^3+5\mathbf{q}+\mathbf{2}$  &\cr
height2pt&\omit&&\omit&&\omit&&\omit& \cr 
& $(4,2,0)$ && $-\mathbf{q}^4+\mathbf{q}^2-\mathbf{q}-\mathbf 1$ && 
$(5,2,0)$  && $2\mathbf{q}^5-\mathbf{q}^4-4\mathbf{q}^3+2\mathbf{q}^2+4\mathbf{q}+\mathbf{1}$  &\cr
height2pt&\omit&&\omit&&\omit&&\omit& \cr 
& $(4,1,1)$ && $-\mathbf{q}^4+\mathbf{q}^2-\mathbf{q}$ && 
$(5,1,1)$ && $\mathbf{q}^5-\mathbf{q}^4-2\mathbf{q}^3+4\mathbf{q}^2+6\mathbf{q}+\mathbf{2}$  &\cr
height2pt&\omit&&\omit&&\omit&&\omit& \cr 
&$(3,3,0)$ && $\mathbf{q}^3+\mathbf{q}^2$ && 
$(4,3,0)$ && $\mathbf{q}^5-2\mathbf{q}^4-3\mathbf{q}^3+2\mathbf{q}^2+3\mathbf{q}+\mathbf{1}$  &\cr
height2pt&\omit&&\omit&&\omit&&\omit& \cr 
&$(3,2,1)$ && $\mathbf{q}^6-\mathbf{q}^4+\mathbf{q}^3+\mathbf{q}^2-\mathbf{q}-\mathbf 1$ && 
$(4,2,1)$ && $-\mathbf{q}^7+2\mathbf{q}^5-2\mathbf{q}^4-3\mathbf{q}^3+3\mathbf{q}^2+3\mathbf{q}$  &\cr
height2pt&\omit&&\omit&&\omit&&\omit& \cr 
&$(2,2,2)$ && $\mathbf{1}$ && 
$(3,3,1)$ && $\mathbf{q}^8+\mathbf{q}^5-2\mathbf{q}^4-2\mathbf{q}^3+2\mathbf{q}^2+3\mathbf{q}+\mathbf{1}$  &\cr
height2pt&\omit&&\omit&&\omit&&\omit& \cr 
&&&  &&  
$(3,2,2)$ && $-\mathbf{q}^7+\mathbf{q}^6+2\mathbf{q}^5-\mathbf{q}^4-\mathbf{q}^3$  &\cr
height2pt&\omit&&\omit&&\omit&&\omit& \cr 
} \hrule}}
\end{thm}

\section{Preliminaries} \label{sec-preliminaries} 
\begin{ntn} We will use the notation $\lambda=[1^{\lambda_1},\ldots,\nu^{\lambda_{\nu}}]$ for the partition where $i$ appears $\lambda_i$ times, of the positive integer $\abs{\lambda}=\sum_{i=1}^{\nu} i \cdot \lambda_i$. The number $\abs{\lambda}$ will be called the weight of the partition $\lambda$ and we will consider the empty set to be the only partition of weight $0$.
\end{ntn}

\begin{dfn} \label{dfn-zee} For a partition $\lambda=[1^{\lambda_1},\ldots,\nu^{\lambda_{\nu}}]$ put 
$\zee{\lambda}:=\prod_{i=1}^{\nu}\lambda_i! \cdot i^{\lambda_i}. $
\end{dfn}

\begin{ntn} Let us fix a finite field $k$ with $q=p^r$ elements. We denote by $k_m$ a degree $m$ extension of $k$ and by $\bar k$ an algebraic closure of $k$.
\end{ntn}

\begin{ntn} \label{ntn-F} Let $X$ be a scheme defined over $k$. The \emph{absolute Frobenius}, denoted by $F'$, is the morphism from $X$ to itself which is the identity on the underlying topological space and the $q$:th power morphism on the sheafs of rings. The morphism $F'\times id$, from  $X_{\bar k} :=X \otimes_k \bar k$ to itself, is called the \emph{$\bar k$-linear Frobenius} morphism. It is indeed $\bar k$-linear, and if $X_{\bar k}$ is a non-singular irreducible projective curve it is a morphism of degree $p$. Let $\phi$ be the automorphism of $\bar k$ that sends $a \in \bar k$ to $a^q \in \bar k$. The \emph{geometric Frobenius} will be denoted by $F$ and it is the morphism from $X_{\bar k}$ to itself that is defined by $id \times \phi^{-1}$. 
\end{ntn} 

\begin{con} 
The word \emph{``Frobenius''} will refer to the geometric Frobenius. 
\end{con} 

\begin{rmk} Observe that even though schemes often will be defined over finite fields, we always consider them over an algebraic closure of the finite field. By remembering the Frobenius morphism we can still talk about its structure over the finite field.
\end{rmk}

Since $F$ is an automorphism, it induces an action on the subschemes of $X_{\bar k}$. 

\begin{dfn} \label{dfn-ord-lambda} Let $X$ be a scheme defined over $k$. An $n$-tuple $(p_1, \ldots, p_n)$ of distinct subschemes of $X_{\bar k}$ is called a \emph{conjugate $n$-tuple} if $Fp_i=p_{i+1}$ for $1 \leq i \leq n-1$ and $Fp_n=p_{1}$. A $\abs{\lambda}$-tuple  $(p_1, \ldots, p_{\abs{\lambda}})$ of distinct subschemes of $X$ is called a \emph{$\lambda$-tuple} if it consists of $\lambda_1$ conjugate $1$-tuples, followed by $\lambda_2$ conjugate $2$-tuples etc.
\end{dfn} 

\begin{dfn} \label{dfn-unord-mu}
For a partition $\mu$, a set $S$ of $\abs{\mu}$ distinct subschemes of $X_{\bar k}$ is called a \emph{$\mu$-set} if there is a $\mu$-tuple $(p_1,\dots,p_{\abs{\mu}})$ such that $S=\{p_1,\dots,p_{\abs{\mu}}\}$. 
If $n$ is a natural number, we will often write $n$-set instead of $[n^1]$-set. Note that an $n$-set is the same as a rational point of degree $n$.
\end{dfn}

\begin{ntn} If not stated otherwise, Greek letters will denote partitions and Roman letters will denote natural numbers.
\end{ntn}

\begin{rmk} \label{rmk-subschemes} 
Let $S$ be a set of subschemes, of a scheme $X_{\bar k}$ defined over $k$, and say that $S$ is fixed by Frobenius, i.e. $S=\{Fs:s \in S\}$. Then $S$ is a $\mu$-set for a unique partition $\mu$ of the number of elements of $S$. 
\end{rmk}

\begin{dfn} \label{dfn-numb-lambda} For any scheme $X$ defined over $k$, let $X(\lambda)$ denote the set of $\lambda$-tuples of points of $X$ whenever $\abs{\lambda} >0$. For $\lambda=\emptyset$ we put $X(\lambda):=\{\emptyset\}$.
\end{dfn}

\begin{lem} \label{lem-Xlambda} Let $\mu$ denote the M\"obius function. For any scheme $X$ defined over $k$, the number of $\lambda$-tuples of points of $X$ is equal to 
\begin{equation} \label{eq-Xlambda}
\abs{X(\lambda)} = \prod_{i=1}^{\nu}\prod_{j=0}^{\lambda_i-1} \Biggl(\biggl(\sum_{d | i} \mu\biggl(\frac{i}{d}\biggr)\cdot\abs{X(k_d)}\biggr) - i \cdot j \Biggr). 
\end{equation}
\end{lem}

\section{Equivariant counts of quartic curves} \label{sec-quartics} 
Let us denote by $\Qqt{n}$ the coarse moduli space of $\Qq{n}$ extended to $\bar{k}$. By an $\s_n$-equivariant count of the number of points defined over $k$ of $\Qqt{n}$ we will mean a count, for each $\sigma \in \s_n$, of the number of fixed points of $F \cdot \sigma$ acting on $\Qqt{n}$. The numbers $\abs{\Qqt{n}^{F \cdot \sigma}}$ only depend upon the cycle type $c(\sigma)$ of the permutation $\sigma$, which is a partition we will often denote by $\lambda$. An $\s_n$-equivariant count of points is called \emph{polynomial} if, for each $\sigma \in \s_n$, there exists a polynomial $f_{\sigma}$ such that $\abs{\Qqt{n}^{F \cdot \sigma}}=f_{\sigma}(q)$ for all but a finite number of characteristics of $k$. 

Define $\mathcal{S}_{\sigma}$ to be the category of non-hyperelliptic genus $3$ curves defined over $k$ together with marked points $(p_1,\ldots,p_n)$ on $C$ defined over $\bar{k}$, such that $(F \cdot \sigma) (p_i)=p_i$ for all $i$. The morphisms of $\mathcal{S}_{\sigma}$ are isomorphisms of the curves that respect the marked tuples of points. 

The points of $\Qqt{n}$ are isomorphism classes of $n$-pointed non-hyperelliptic genus $3$ curves defined over $\bar k$. If the pointed curve $C$ is a representative of a point in $\Qqt{n}^{F \cdot \sigma}$ there is an isomorphism from $C$ to the pointed curve $(F \cdot \sigma)C$. This isomorphism gives a way of descending to an object of $\mathcal{S}_{\sigma}$ (see Lemma 10.7.5 in \cite{Katz-Sarnak}). The number of $\bar k$-isomorphism classes of the category $\mathcal{S}_{\sigma}$ is therefore equal to $\abs{\Qqt{n}^{F \cdot \sigma}}$.

Fix an object $(C,p_1,\ldots,p_n)$ in $\mathcal{S}_{\sigma}$. The sum, over all $k$-isomorphism classes of $n$-pointed curves $(D,q_1,\ldots,q_n)$ that are $\bar{k}$-isomorphic to $(C,p_1,\ldots,p_n)$, of the reciprocal of the number of $k$-automorphisms of $(D,q_1,\ldots,q_n)$, is equal to $1$ (see \cite{Geer} or Lemma 10.7.5 in \cite{Katz-Sarnak}). This enables us to go from $\bar k$-isomorphism classes to $k$-isomorphism classes:
$$\abs{\Qqt{n}^{F \cdot \sigma}} = \sum_{[Y] \in \mathcal{S}_{\sigma}/\cong_{\bar{k}}} 1 =  \sum_{[Y] \in \mathcal{S}_{\sigma}/\cong_{\bar{k}}} \sum_{\substack{[X] \in \mathcal{S}_{\sigma}/\cong_{k} \\ X\cong_{\bar{k}}Y}} \frac{1}{\abs{\Aut_{k}(X)}}  = \sum_{[X] \in \mathcal{S}_{\sigma}/\cong_{k}} \frac{1}{\abs{\Aut_{k}(X)}}.$$

Fix a curve $C$ in $\Qq{}(k)$ and let $X_1, \ldots, X_m$ be representatives of the distinct $k$-isomorphism classes of the subcategory of $\mathcal{S}_{\sigma}$ of elements $(D,q_1,\ldots,q_n)$ such that $D \cong_k C$. For each $X_i$ we can act with $\Aut_k(C)$. This gives an orbit lying in $\mathcal{S}_{\sigma}$ and the stabilizer of $X_i$ is equal to $\Aut_k(X_i)$. Together the orbits of $X_1, \ldots, X_m$ will contain $\abs{C\bigl(c(\sigma)\bigr)}$ elements (recall Definitions \ref{dfn-ord-lambda} and \ref{dfn-numb-lambda}) and hence we obtain
\begin{equation*} 
\abs{\Qqt{n}^{F \cdot \sigma}} = \sum_{[X] \in \mathcal{S}_{\sigma}/\cong_{k}} \frac{1}{\abs{\Aut_{k}(X)}} = \sum_{[C] \in \Qq{}(k)/\cong_k} \frac{\abs{C\bigl(c(\sigma)\bigr)}}{\abs{\Aut_k(C)}}.  
\end{equation*}

Using the canonical divisor we can embed the curves of $\Qq{}(k)$ as degree $4$ curves in the projective plane. After a choice of coordinates over $k$, each plane curve of degree $4$ over $k$ can be represented by a point in the $\Pb^{14}(k)$ of its coefficients. Let $Q(k)$ denote the subset of $\Pb^{14}(k)$ consisting of non-singular curves.

The $k$-automorphisms of the curves in $Q(k)$ are precisely given by the elements of $\mathrm{PGL}_{3}(k)$. This gives an action of $\mathrm{PGL}_{3}(k)$ on $Q(k)$ and the stabilizer of $C$ in $Q(k)$ is thus equal to $\Aut_k(C)$. We can then conclude that $\abs{\Qqt{n}^{F \cdot \sigma}}$ is equal to 
$$\sum_{[C] \in \Qq{}(k)/\cong_k} \frac{\abs{C\bigl(c(\sigma)\bigr)}}{\abs{\Aut_k(C)}} = \sum_{[C] \in Q(k)/\mathrm{PGL}_3(k)} \frac{\abs{C\bigl(c(\sigma)\bigr)}}{\abs{\mathrm{Stab}(C)}}  =  \frac{1}{\PGLk{3}} \sum_{C \in Q(k)} \abs{C(c(\sigma))}.$$
This reduces the question of computing $\abs{\Qqt{n}^{F \cdot \sigma}}$ to computing the number of $\lambda$-tuples on each curve of $Q(k)$. Or equivalently, to computing, for each choice of a $\lambda$-tuple of points in the plane, the number of curves in $Q(k)$ that pass through it.

\section{The sieve principle} \label{sec-sieve}  
\begin{dfn} \label{dfn-Lp}
Let us from now on identify the space of plane degree $4$ curves defined over $k$ with $\Pb^{14}(k)$. For every $\lambda$-tuple $P$ of points in the projective plane, we define $L_P$ to be the linear subspace of $\Pb^{14}(k)$ of curves that contain $P$. 
\end{dfn}

Fix a $\lambda$-tuple $P$ of points in the plane. We would like to find how many non-singular quartic curves over $k$ there are that contain $P$. We will begin this computation with the use of the sieve principle, starting with the space $L_P$ and then successively adding and removing loci of singular curves in such a way that in each step, more and more singular curves will have been removed precisely once.

\begin{dfn}\label{dfn-Lps} 
For any $\mu$-set $S$ of points in the plane, we denote by $L_{P,S}$ the linear subspace of $L_P$ of curves that have singularities at the points of $S$.
\end{dfn}

The locus of singular curves in $L_P$ is the union of all linear spaces $L_{P,S}$ for every $m\geq 1$ and $m$-set $S$. We want to use the sieve principle to compute the number of elements of this union. That is, sum the numbers $(-1)^{i+1}  \abs{L_{P,S_1} \cap \ldots \cap L_{P,S_i}}$, for each $i \geq 1$ and for each unordered choice of distinct sets $S_1, \ldots, S_i$ where $S_j$ is an $m_j$-set. If this procedure terminates, every singular curve in $L_P$ will have been counted exactly once. Hence, taking $\abs{L_P}$ minus the resulting number gives the number of non-singular curves in $L_P$.

Note that, if $S=S_1 \cup S_2$ for a $\lambda$-set $S_1$ and a $\mu$-set $S_2$, then $L_{P,S_1} \cap L_{P,S_2} = L_{P,S}$. Thus, to be able to use the sieve principle in this way we need only find the dimensions of all linear spaces of the form $L_{P,S}$.

Unfortunately, determining the dimension of all linear subspaces of the form $L_{P,S}$  is not always easy. Moreover, in order to apply the sieve principle we need to know that there is a number $M$ such that $L_{P,S}$ is empty as soon as $S$ consists of more than $M$ points. Since there are curves with infinitely many singularities, namely the nonreduced ones, such a number $M$ does not exist.
Instead, we choose a number $M$ and define a modified sieve principle as follows.

\begin{dfn}
Let $P$ be a $\lambda$-tuple of points in the plane, and $M$ a positive integer. We define the \emph{modified sieve principle} as the computation of the sum of $(-1)^{i+1} \abs{L_{P,S_1} \cap \ldots \cap L_{P,S_i}}$ for each $i \geq 1$ and for each unordered choice of distinct sets $S_1, \ldots, S_i$ where $S_j$ is an $m_j$-set such that $\sum_{j=1}^i m_j \leq M$. Denote the resulting number by $s_{M,P}$.
\end{dfn}

We see that if we subtract $s_{M,P}$ from $\abs{L_P}$, then all curves with at most $M$ singularities will have been removed from $L_P$ exactly once. 

\begin{dfn} \label{dfn-linsub} Define $\linsub{M}:=\sum_{P}(\abs{L_P}-s_{M,P})$ where the sum runs over all $\lambda$-tuples $P$ of points in the plane.
\end{dfn}

Thus, if we have computed $\linsub{M}$ we need to amend for the curves with more than $M$ singularities in order to obtain the sum over all $\lambda$-tuples $P$ of the number of non-singular curves in $L_P$. 

\begin{dfn} \label{dfn-t} For any partition $\mu$ of weight $m>M$ and any partition $\lambda$, define $t_{\lambda,\mu}$ to be the sum over all choices of $\lambda$-tuples $P$ and $\mu$-sets $S$ of the number of curves that contain $P$ and that have singularities at the points of $S$ and nowhere else. 

In this definition we also allow $\mu$ to be an infinite partition, that is, of the kind $\mu=[1^{\mu_1},\ldots,\nu^{\mu_{\nu}},\ldots]$.
\end{dfn}

\begin{rmk} 
If $\lambda$ is fixed, the $t_{\lambda,\mu}$'s are nonzero only for a finite number of choices of $\mu$. 
\end{rmk}

For every curve with more than $M$ singularities we can easily find how many times it would be removed or added when applying the modified sieve principle to compute $s_{M,P}$. We see that the multiplicity with which a curve with more than $M$ singularities is counted in the computation of $s_{M,P}$ is the same for all curves which are singular exactly at a $\mu$-set. Therefore, after adding a suitable multiple of $t_{\lambda,\mu}$ to $\linsub{M}$ for each $\mu$ with $\abs{\mu} > M$, every singular curve will have been removed exactly once from each space $L_P$ and hence we will have computed $\abs{\Qqt{n}^{F \cdot \sigma}}$. 

\subsection{The choice of $M$} \label{sec-M1}
There are two parts to the method to make $\s_n$-equivariant counts of $\Qq{n}$ that was introduced above. For a partition $\lambda$, we should on the one hand compute $\linsub{M}$ and on the other $t_{\lambda,\mu}$ for all $\mu$ with $\abs{\mu} > M$. The choice of $M$ is done by weighing the difficulty of these two parts. 

We will choose $M=1$ and make $\s_n$-equivariant counts of the number of points of $\Qq{n}$ for $n \leq 7$. Thus, we need to find the dimensions of $L_{P,S}$ for all $\lambda$-tuples $P$ and $\mu$-sets $S$, where $\abs{\lambda} \leq 7$ and $\abs{\mu} \leq 1$. This is carried out in Section \ref{sec-linsub}. 
We also need to compute $t_{\lambda,\mu}$ when $\abs{\lambda} \leq 7$ and $\abs{\mu} \geq 2$, which is done in Sections \ref{sec-library} to~\ref{sec-two}. After having computed $\linsub{1}$ we find that the singular curves have been removed once for each singularity over $k$. This shows that 
\begin{multline*} 
\abs{\Qqt{n}^{F \cdot \sigma}} = \frac{1}{\PGLk{3}}(\linsub{1} +  \sum_{\abs{\mu} > 1 } (\mu_1-1) \cdot t_{\lambda,\mu}) =\\= \frac{1}{\PGLk{3}}\Bigl(\sum_{P \in \Pb^2(\lambda)} \bigl(\abs{L_{P}} - \sum_{S \in \Pb^2([1^1])} \abs{L_{P,S}} \bigr) +  \sum_{\abs{\mu} > 1 } (\mu_1-1) \cdot t_{\lambda,\mu} \Bigr). 
\end{multline*}

\section{The library of singular degree $4$ curves} \label{sec-library} 
\begin{con} 
An \emph{irreducible} curve in $\Pb^2_{\bar k}$ of degree $d=1$ (resp. $d=2,3,4$) will be called a line (resp. conic, cubic, quartic).
\end{con}

\begin{dfn} \label{dfn-types} 
We define the \emph{type} of a degree $4$ curve over $k$ to be given by the following information: 
\begin{itemize}
\item[$\star$] degrees and multiplicities of the irreducible components; 
\item[$\star$] over which fields the components are defined; 
\item[$\star$] the number of singularities of each irreducible component;
\item[$\star$] the delta invariants $\delta$ of the singularities of each irreducible component;
\item[$\star$] over which fields the singularities of each irreducible component are defined;
\item[$\star$] the number of points in the inverse image, under the normalization morphism, of each singularity of an irreducible component; 
\item[$\star$] over which fields the points of the inverse images are defined;
\item[$\star$] in how many points every set of irreducible components intersect;
\item[$\star$] over which fields the intersection points are
defined. 
\end{itemize}
\end{dfn}

In Sections \ref{sec-inf} to \ref{sec-two} we will go through the different types of quartic curves with a $\mu$-set of singularities, when $\abs{\mu} \geq 2$, and find sufficient information to compute their contribution to $t_{\lambda,\mu}$ when $\abs{\lambda} \leq 7$. 

It follows from Definition \ref{dfn-types} that knowing the type of a curve together with the number of points over $k_n$, for all $n$, of the normalization of each of its components, is sufficient to compute its contribution to $t_{\lambda,\mu}$ for \emph{any} $\lambda$. A non-singular curve over $k$ of genus $0$ is isomorphic to $\Pbk^1$ and $\abs{\Pbk^1(k_n)}=q^n+1$ for all $n$. Moreover, knowing the intersection points of a conjugate $m$-tuple of non-singular genus 0 curves is sufficient to conclude its number of points over any extension of $k$. We conclude from this that for all types consisting of curves with irreducible components whose normalization have genus $0$, we only need to find how many curves there are of the type considered. 

There are also types of curves with a component whose normalization has genus~$1$. First we have the non-singular cubic curves together with a line. We will find their contribution by formulating these counts in terms of $\s_n$-equivariant counts of points of $\Mm{1}{n}$ for $n \geq 1$. In turn, the $\s_n$-equivariant counts of points of $\Mm{1}{n}$ will be found by computing the number of pointed non-singular cubics with the same method as we are using to count the pointed non-singular quartics. This is done in Section~\ref{sec-cubics} for $n \leq 10$. Alternatively, the information on $\Mm{1}{n}$, for any $n$, can be found using the results of \cite{G-res}.

Then we have the types consisting of quartics with two singularities with $\delta=1$. We will find a set of morphisms that contains the normalization morphisms of these curves, and the number of elements of this set can be computed in terms of $\s_n$-equivariant counts of $\Mm{1}{n}$. But to count the morphisms of this set that have an image of the wrong kind, which will be the morphisms of degree $2$, we will need somewhat more detailed information on $\Mm{1}{n}$. To find this information we will employ methods and results found in \cite{Jonas1}.

\subsection{Conventions for tables} \label{sec-conv} 
There will be tables below, presenting the results on the number of degree $4$ curves of each type where the irreducible components are of degree $1$ or $2$. In these tables we will use the following conventions: 
\begin{itemize}
\item[$\star$] components are defined over the ground field $k$, if nothing else is stated;
\item[$\star$] if we have $n$ intersection (abreviated: int.) points we will write out the partition they induce on $n$. 
\end{itemize}
The column marked ``$\#$'' contains the number of curves, of each type, divided by $\PGLk{3}$. Actually all results on numbers of curves presented below will be after dividing by $\PGLk{3}$.

\section{Infinitely many singularities} \label{sec-inf} 
The curves with infinitely many singularities are precisely the non-reduced ones, and all such will have irreducible components of degree at most $2$. 

\begin{exa}\label{exa-conic} 
Let us count the number of curves of the type, double conics. There are $\abs{\Pb^5(k)}$ degree $2$ curves over $k$ and the reducible ones are of the following types: pairs of distinct lines over $k$, conjugate pairs of lines, or double lines over $k$. Hence we get, 
\begin{equation*}\bigl(\frac{q^6-1}{q-1}-\binom{q^2+q+1}{2}-\frac{1}{2}(q^4-q)-(q^2+q+1)\bigr) \frac{1}{\PGLk{3}}= \frac{1}{q(q+1)(q-1)}.
\end{equation*}
\end{exa}

\begin{exa}
Consider the curves of the type, a conic and a double line, both defined over $k$, that intersect in a conjugate pair. Choose a conjugate pair of points $(p_1,p_2)$ in the plane and let them define the double line, which necessarily will be over $k$. Then we should pick a conic over $k$ in the $\abs{\Pb^3(k)}$ degree $2$ curves over $k$ that contain $(p_1,p_2)$. The reducible curves are pairs of lines of two kinds. Either take any line over $k$ together with the line through $p_1$, $p_2$. Or take any point $q$ over $k$ outside the line through $p_1$, $p_2$ and let the lines be the ones through $p_1$, $q$ and $p_2$, $q$. Putting the pieces together we get 
$$\frac{1}{2}(q^4-q) (\frac{q^4-1}{q-1}-(q^2+q+1)-q^2) \frac{1}{\PGLk{3}}=\frac{1}{2(q+1)}.$$
\end{exa}

\begin{table}[htbp] \caption{Infinitely many singularities} \label{tab-inf}
\centerline{
\vbox{
\offinterlineskip
\hrule
\halign{&\vrule#& \quod #\hfil \strut &\vrule#& \quod \hfil#\hfil \quod \strut  \cr
height2pt&\omit&&\omit& \cr
& \bf{Description} && \bf{\#} & \cr 
height2pt&\omit&&\omit&\cr
\noalign{\hrule}
height2pt&\omit&&\omit&\cr
&\emph{\small Quadruple line} &&$\scriptstyle \frac{1}{q^3(q+1)(q-1)^2}$&\cr
height2pt&\omit&&\omit&\cr
&\emph{\small Line + triple line, int. in $\scriptstyle [1^1]$} &&$\scriptstyle \frac{1}{q^2(q-1)^2}$ &\cr
height2pt&\omit&&\omit&\cr
&\emph{\small $\lambda$-set of double lines, int. in $\scriptstyle [1^1]$} && &\cr
height2pt&\omit&&\omit&\cr
&$\scriptstyle \star \,\, \lambda=[1^2],[2^1]$ &&$\scriptstyle \frac{1}{2q^2(q-1)^2},\frac{1}{2q^2(q-1)(q+1)}$&\cr
height2pt&\omit&&\omit&\cr
&\emph{\small Double conic} &&$\scriptstyle \frac{1}{q(q+1)(q-1)}$ &\cr
height2pt&\omit&&\omit&\cr
&\emph{\small $\lambda$-set of lines + double line, int. in $\mu$-set} && &\cr
height2pt&\omit&&\omit&\cr
&$\scriptstyle \star \,\, (\lambda,\mu)=([1^2],[1^1]), ([2^1],[1^1]), ([1^2],[1^3]), ([2^1],[1^1,2^1])$ &&$\scriptstyle \frac{1}{2q^2(q-1)},\frac{1}{2q^2(q-1)},\frac{1}{2(q-1)^2},\frac{1}{2(q+1)(q-1)}$&\cr
height2pt&\omit&&\omit&\cr
&\emph{\small Conic + double line, int. in $\scriptstyle [1^1], [1^2], [2^1]$} &&$ \scriptstyle \frac{1}{q(q-1)}, \frac{1}{2(q-1)}, \frac{1}{2(q+1)}$ &\cr
height2pt&\omit&&\omit&\cr
} \hrule}}
\end{table}

\section{Six singularities} \label{sec-six} The only possibility to get six singularities is to have four lines in general position, that is, where no three of them intersect in one point. 
The type is determined by choosing the fields over which the lines are defined, because this, in turn, determines over which fields the intersection points are defined.

\begin{table}[htbp] \caption{Six singularities } \label{tab-six}
\centerline{
\vbox{
\offinterlineskip
\hrule
\halign{&\vrule#& \quod #\hfil \strut &\vrule#& \quod \hfil#\hfil \quod \strut  \cr
height2pt&\omit&&\omit& \cr
& \bf{Description} && \bf{\#} & \cr
height2pt&\omit&&\omit&\cr
\noalign{\hrule}
height2pt&\omit&&\omit&\cr
&\emph{\small $\lambda$-set of lines, int. in $\mu$-set} && &\cr 
height2pt&\omit&&\omit&\cr
&$\scriptstyle \star \,\, (\lambda,\mu)=([1^4],[1^6]), ([1^2,2^1],[1^2,2^2]), ([2^2],[1^2,2^2]), ([1^1,3^1],[3^2]), ([4^1],[2^1,4^1])$ &&$\scriptstyle \frac{1}{24},\frac{1}{4},\frac{1}{8},\frac{1}{3},\frac{1}{4}$&\cr
height2pt&\omit&&\omit&\cr
} \hrule}}
\end{table}

\section{Five singularities} \label{sec-five} The degree $4$ curves with five singularities all consist of a conic and two lines that are in general position. The type of such a curve follows from the choice of the fields over which the lines and the intersection points are defined.

\begin{exa} \label{exa-51} Let us compute how many choices there are of a conic and two lines intersecting in five points over $k$. Due to the duality of lines and points in the plane, we can read off from Table~\ref{tab-six} that there are $1/24 \cdot \PGLk{3}$ possibilities for four general points over $k$. Let us choose such a four-tuple $P$, and then one of the three different pairs of lines over $k$ that contain $P$. The two chosen lines will then intersect in a point defined over $k$, disjoint from $P$. Finally, we will choose a conic in the $\Pb^1$ of degree $2$ curves that contain $P$. We already know that there are three reducible curves over $k$ in this $\Pb^1$ and hence we get
$$\frac{1}{24}\PGLk{3} \cdot 3 (q+1-3) \frac{1}{\PGLk{3}} = \frac{1}{8}(q-2).$$
\end{exa}

\begin{table}[htbp] \caption{Five singularities} \label{tab-five}
\centerline{
\vbox{
\offinterlineskip
\hrule
\halign{&\vrule#& \quod #\hfil \strut &\vrule#& \quod \hfil#\hfil \quod \strut  \cr
height2pt&\omit&&\omit& \cr
& \bf{Description} && \bf{\#} & \cr
height2pt&\omit&&\omit&\cr
\noalign{\hrule}
height2pt&\omit&&\omit&\cr
&\emph{\small Conic + $\lambda$-set of lines, int. in $\mu$-set} && &\cr 
height2pt&\omit&&\omit&\cr
&$\scriptstyle \star \,\, (\lambda,\mu)=([1^2],[1^5]), ([1^2],[1^3,2^1]), ([1^2],[1^1,2^2])$ &&$\scriptstyle \frac{1}{8}(q-2),\frac{1}{4}q, \frac{1}{8}(q-2)$&\cr
height2pt&\omit&&\omit&\cr
&$\scriptstyle \star \,\, (\lambda,\mu)=([2^1],[1^1,2^2]), ([2^1],[1^1,4^1])$ &&$\scriptstyle \frac{1}{4}(q-2),\frac{1}{4}q$&\cr
height2pt&\omit&&\omit&\cr
} \hrule}}
\end{table}

\section{Four singularities} \label{sec-four}
First we have 
types of curves with the same kinds of components as in the types with five and six singularities, but where the components no longer are in general position. Then there is the new case of two conics intersecting transversally.

\begin{exa} Let us compute the number of curves of type, a conjugate pair of conics intersecting in a conjugate quadruple. We begin by choosing a conjugate quadruple of points $(p_1,p_2,p_3,p_4)$, no three of which lie on a line. Due to the duality of lines and points in the plane, we can read off from Table~\ref{tab-six} that there are $1/4 \cdot \PGLk{3}$ such. There is a $\Pb^1$ of degree $2$ curves that contain $p_1,p_2,p_3,p_4$, and we want to count the conjugate pairs of such that are conics. There is only one conjugate pair of reducible degree $2$ curves in this $\Pb^1$. Namely the pair of lines through $p_1$, $p_2$ and $p_3$, $p_4$, together with the pair of lines through $p_1$, $p_4$ and $p_2$, $p_3$. We conclude that the number we seek is
$$\frac{1}{4} \PGLk{3} (\frac{1}{2}(q^2-q)-1) \frac{1}{\PGLk{3}} = \frac{1}{8}(q^2-q-2).$$
\end{exa}

\begin{table}[htbp] \caption{Four singularities}\label{tab-four}
\centerline{
\vbox{
\offinterlineskip
\hrule
\halign{&\vrule#& \quod #\hfil \strut &\vrule#& \quod \hfil#\hfil \quod \strut  \cr
height2pt&\omit&&\omit& \cr
& \bf{Description} && \bf{\#} & \cr
height2pt&\omit&&\omit&\cr
\noalign{\hrule}
height2pt&\omit&&\omit&\cr 
&\emph{\small $\lambda$-set of lines, int. in $\mu$-set} && &\cr
height2pt&\omit&&\omit&\cr
&$\scriptstyle \star \,\, (\lambda,\mu)=([1^4],[1^4]),([1^2,2^1],[1^2,2^1]),([1^1,3^1],[1^1,3^1])$ &&$ \scriptstyle \frac{1}{6(q-1)},\frac{1}{2(q-1)},\frac{1}{3(q-1)}$ &\cr
height2pt&\omit&&\omit&\cr
&\emph{\small Conic + tgt line + transv. line, int. in $\lambda$-set} && &\cr 
height2pt&\omit&&\omit&\cr
&$\scriptstyle \star \,\,\lambda=[1^4],[1^2,2^1]$ &&$\frac{1}{2}$,$\frac{1}{2}$&\cr
height2pt&\omit&&\omit&\cr
&\emph{\small $\lambda$-set of conics, int. in $\mu$-set} && &\cr
height2pt&\omit&&\omit&\cr
&$\scriptstyle \star \,\, (\lambda,\mu)=([1^2],[1^4]),([1^2],[1^2,2^1]),([1^2],[2^2])$ &&$ \scriptstyle \frac{1}{48}(q-2)(q-3), \frac{1}{8}q(q-1), \frac{1}{16}(q-2)(q-3)$ &\cr
height2pt&\omit&&\omit&\cr
&$\scriptstyle \star \,\, (\lambda,\mu)=([1^2],[1^1,3^1]),([1^2],[4^1])$ &&$ \scriptstyle \frac{1}{6}q(q+1), \frac{1}{8}q(q-1)$ &\cr
height2pt&\omit&&\omit&\cr
&$\scriptstyle \star \,\, (\lambda,\mu)=([2^1],[1^4]),([2^1],[1^2,2^1]),([2^1],[2^2])$ &&$ \scriptstyle \frac{1}{48}(q^2-q),\frac{1}{8}(q^2-q-2),\frac{1}{16}(q^2-q)$&\cr
height2pt&\omit&&\omit&\cr
&$\scriptstyle \star \,\, (\lambda,\mu)=([2^1],[1^1,3^1]),([2^1],[4^1])$ &&$ \scriptstyle \frac{1}{6}(q^2-q),\frac{1}{8}(q^2-q-2)$ &\cr
height2pt&\omit&&\omit&\cr
} \hrule}}
\end{table}

\subsection{Singular cubic with transversal line} \label{sec-cubtrans}
A cubic curve can have at most one singularity, it must have $\delta =1$, and if we resolve the singularity we get a non-singular curve of genus $0$. Let us define the necessary distinctions of singularities with $\delta =1$ and $\delta =2$ that we will use throughout the article. 

\begin{dfn} \label{dfn-sings} 
Say that we have an $i$-set of singularities with $\delta = 1$ on a curve over $k$. The inverse image under the normalization morphism of the singularities will be a $\lambda$-set where $\lambda$ is either $[i^1]$, $[i^2]$ or $[2i]$. The singularities will be called \emph{cusps}, \emph{split nodes}, \emph{non-split nodes} and be denoted $c^i$, $n_1^i$ and $n_2^i$ respectively. If $\delta=2$ then we make the same distinction, but the singularities will be called \emph{taccusps}, \emph{split tacnodes} or \emph{non-split tacnodes} and be denoted $tc^i$, $tn_1^i$, $tn_2^i$ respectively. Moreover, we shall suppress $i$ in the notation when $i=1$. 
\end{dfn}

From this characterisation we see that a cuspidal cubic curve has the same number of points over every extension of $k$ as $\Pbk^1$, the split nodal as many as $\Pbk^1$ minus one point over $k$, and finally the non-split nodal as many as $\Pbk^1$ plus one point over $k$ and minus a conjugate pair of points.

What we are left to find is the number of curves of each type that consist of a singular cubic with a transversal line. 

\begin{lem} \label{lem-pgl2} A singular cubic with fixed normalization has $\PGLk{2}$ distinct normalization morphisms over $k$.
\end{lem}
\begin{pf}
The normalization morphism is unique up to unique isomorphism. 
\end{pf}
 
As we have seen, the inverse image of the singularity under the normalization morphism of the cubic, is a $\lambda$-set $P$ where $\abs{\lambda}=1$ or $2$. 
\begin{rmk} \label{rmk-abuse} With slight abuse of notation, we will write $P=\{p_1,p_2\}$ also in the cuspidal case, even though in this case $p_1$ and $p_2$ are equal.
\end{rmk}
The inverse image of the intersection points of the transversal line with the singular cubic gives a $\mu$-set $Q$ where $\abs{\mu}=3$. Note that the type of the curve is determined by the partitions $\lambda$ and $\mu$. 

We will now construct the normalization morphisms using such choices of points on $\Pbk^1$. Fix a partition $\lambda$ of weight $1$ or $2$, and a partition $\mu$ of weight $3$. Then choose a $\lambda$-set $P=\{p_1,p_2\}$ and a $\mu$-set $Q=\{q_1,q_2,q_3\}$ of points on $\Pbk^1$ such that $P$ and $Q$ are disjoint. Let $\mathfrak{d}$ be the linear system inside the complete linear system $|q_1+q_2+q_3|$ spanned by the divisor $q_1+q_2+q_3$ and all divisors of the form $p_1+p_2+r$, where $r$ may be any point on $\Pbk^1$. The divisor $q_1+q_2+q_3$ lies outside the line defined by $p_1+p_2+r$ and therefore $\mathfrak{d}$ has dimension $2$. Since there clearly are divisors in $\mathfrak{d}$ that have no point in common, $\mathfrak{d}$ is base point free.

\begin{lem} \label{lem-pgl3} From a base point free linear system $\mathfrak{d}'$ of dimension $2$ on a curve $X$ defined over $k$, we get $\PGLk{3}$ distinct morphisms to $\Pbk^2$. 
\end{lem}

\begin{lem} \label{clm-cubtrans} The image of any morphism induced by $\mathfrak{d}$ is a singular cubic with $q_1$, $q_2$ and $q_3$ lying on a line. 
\end{lem}
\begin{pf} Fix any morphism $\varphi$ induced by $\mathfrak{d}$. Since the divisors of $\mathfrak{d}$ have degree $3$, the degree of $\varphi$ is either $1$ or $3$. In the latter case, the image would be a line, which is not possible since $\mathfrak{d}$ has dimension $2$. The image of the morphism is therefore a cubic and since the $\lambda$-set $P$ lies on all the lines corresponding to the divisors $p_1+p_2+r$ where $r$ is any point on $\Pbk^1$, it will be a singularity. Hence, $\varphi$ is a normalization morphism and the image of the $\mu$-set $Q$ consists of three distinct points that lie on the line corresponding to the divisor $q_1+q_2+q_3$. \end{pf}

This gives us two ways of counting the same set:
\begin{itemize}
\item[$\star$] for each choice of a $\lambda$-set $P$ and $\mu$-set $Q$ on $\Pbk^1$ such that $P$ and $Q$ are disjoint, there are $\PGLk{3}$ normalization morphisms over $k$;
\item[$\star$] for each choice of a singular cubic with a transversal line, of a type determined by the partitions $\lambda$ and $\mu$, there are $\PGLk{2}$ normalization morphisms.
\end{itemize}

This proves the following lemma.

\begin{lem} \label{lem-cubtrans} With notation as above, the number of singular cubics with a trans-versal line of a type determined by the partitions $\lambda$ and $\mu$, is equal to 
$$\frac{1}{\zee{\lambda} \zee{\mu}} \cdot \abs{\Pbk^1(\lambda)} \cdot \abs{X(\mu)} \cdot\frac{1}{\PGLk{2}} \;\; \text{where} \;\; X=\Pbk^1 \setminus \{a \; \lambda-set \}.$$ 
\end{lem}

\section{Three singularities} \label{sec-three} There are types of curves with three singularities consisting of two conics, or consisting of a conic with two lines, where in both these cases the components will not be in general position. 

\begin{table}[htbp] \caption{Three singularities} \label{tab-three-red}
\centerline{
\vbox{
\offinterlineskip
\hrule
\halign{&\vrule#& \quod #\hfil \strut &\vrule#& \quod \hfil#\hfil \quod \strut  \cr
height2pt&\omit&&\omit& \cr
& \bf{Description} && \bf{\#} & \cr
height2pt&\omit&&\omit&\cr
\noalign{\hrule}
height2pt&\omit&&\omit&\cr
&\emph{\small Conic + $\lambda$-set of tgt lines, int. in $\mu$-set} && &\cr 
height2pt&\omit&&\omit&\cr
&$\scriptstyle \star \,\,(\lambda,\mu)=([1^2],[1^3]),([2^1],[1^1,2^1])$ &&$\scriptstyle \frac{1}{2(q-1)},\frac{1}{2(q+1)} $&\cr
height2pt&\omit&&\omit&\cr
&\emph{\small Conic + $\lambda$-set of lines, int. in $\mu$-set on the conic} && &\cr
height2pt&\omit&&\omit&\cr
&$\scriptstyle \star \,\,(\lambda,\mu)=([1^2],[1^3]),([2^1],[1^1,2^1])$ &&$\scriptstyle \frac{1}{2},\frac{1}{2} $&\cr
height2pt&\omit&&\omit&\cr
&\emph{\small $\lambda$-set of conics, int. in $\mu$-set} && &\cr
height2pt&\omit&&\omit&\cr
&$\scriptstyle \star \,\,(\lambda,\mu)=([1^2],[1^3]),([1^2],[1^1,2^1]),([2^1],[1^3]),([2^1],[1^1,2^1])$ &&$\scriptstyle \frac{1}{4}(q-2),\frac{1}{4}(q-2),\frac{1}{4}q,\frac{1}{4}q$&\cr
height2pt&\omit&&\omit&\cr
} \hrule}}
\end{table}

\subsection{Singular cubic with a line intersecting in two smooth points on the cubic} \label{sec-cubtwo} 
This case follows in the same way as in Section~\ref{sec-cubtrans}, with the difference that $Q=(q_1,q_3)$ is a $[1^2]$-tuple and $q_1=q_2$.

\begin{lem} \label{lem-cubtwo} The number of singular cubics with a line intersecting in two smooth points on the cubic, of a type determined by the partition $\lambda$, is equal to 
$$\frac{1}{\zee{\lambda}} \abs{\Pbk^1(\lambda)} \cdot \abs{X([1^2])} \cdot \frac{1}{\PGLk{2}} \;\; \text{where} \;\; X=\Pbk^1 \setminus \{a \; \lambda-set \}. $$
\end{lem}

\subsection{Non-singular cubic with a transversal line} \label{sec-noncubtrans}
These are curves that have a component of both arithmetic and geometric genus equal to $1$, which has the consequence that the number of points of the curves is not determined by the type. We are therefore forced to use a different approach than previously and it will be similar to the one used by Belorousski in \cite{Pashathesis}.

Fix a partition $\mu$ of weight $3$. Let $C$ be a genus $1$ curve and let $P=\{p_1,p_2,p_3\}$ be a $\mu$-set of points on $C$. 

\begin{lem} \label{clm-noncubtrans} The image of any morphism induced by the complete linear system $|p_1+p_2+p_3|$ is a non-singular cubic with $p_1$, $p_2$ and $p_3$ lying on a line. 
\end{lem}

Thus we have a way of realizing the elements of the moduli space $\Mm{1}{3}$ as non-singular cubics with three points lying on a line. Let us spell out this connection in terms of $\s_n$-equivariant counts of points, as in Section \ref{sec-quartics} for $\Qq{}$ and quartics.
For all $n \geq 1$, let $\Mmt{1}{n}$ be the coarse moduli space of $\Mm{1}{n}$ extended to $\bar k$. Then let $\mathcal{C}_{\mu}(k)$ be the set of non-singular cubics over $k$ together with a $\mu$-set of points lying on a line.

\begin{ntn} For each partition $\lambda$, fix a permutation $\sigma_{\lambda}$ with the property $c(\sigma_{\lambda})=\lambda$. 
\end{ntn}

\begin{dfn} \label{dfn-add-subtr} For any partitions $\lambda$ and $\mu$ we denote by $\lambda+\mu$ the partition $[1^{\lambda_1+\mu_1},\ldots,\nu^{\lambda_{\nu}+\mu_{\nu}}]$. We write $\mu \leq \lambda$ if $\mu_i \leq \lambda_i$ for all $i$. If $\mu \leq \lambda$ we denote by $\lambda-\mu$ the partition $[1^{\lambda_1-\mu_1},\ldots,\nu^{\lambda_{\nu}-\mu_{\nu}}]$.
\end{dfn}

Let $\lambda$ be a partition of weight $n$. Following the arguments of Section \ref{sec-quartics} and using Lemma \ref{clm-noncubtrans} we can prove that 
\begin{equation} \label{eq-noncub} 
\frac{1}{\zee{\mu}} \cdot \abs{\Mmt{1}{3+n}^{F \cdot \sigma_{\mu+\lambda} } } = \frac{1}{\PGLk{3}} \cdot \sum_{C \in \mathcal{C}_{\mu}(k)} \abs{X(\lambda)} \;\; \text{where} \;\; X=C \setminus \{the \; \mu-set \}.
\end{equation} 

\begin{lem} \label{lem-noncubtrans} The number of non-singular cubics with a transversal line, of a type determined by the partition $\mu$, and with a $\lambda$-tuple of points is equal to
\begin{equation} \label{eq-noncubtrans} 
\frac{1}{\zee{\mu}} \cdot \sum_{i=0}^{\abs{\lambda}} \sum_{\substack{\kappa \vdash i \\ \kappa \leq \lambda}}\frac{\zee{\lambda}}{\zee{\kappa} \zee{\lambda-\kappa}} \abs{\Mmt{1}{3+i}^{F \cdot \sigma_{\mu+\kappa}} } \cdot \abs{\Pbk^1(\lambda-\kappa)}.\end{equation}
\end{lem}
\begin{pf} For each $\kappa \leq \lambda$ we can distribute a $\lambda-\kappa$-tuple of points on the line, and a $\kappa$-tuple of points on the cubic but outside the $\mu$-set of intersection points. From equation \eqref{eq-noncub} we can find the choices of points on the cubic in terms of the moduli space of genus $1$ curves. Equation \eqref{eq-noncubtrans} follows when taking into account the number of ways of reordering the $\abs{\lambda}$ points. 
\end{pf}

In Section \ref{sec-cubics} we will make $\s_n$-equivariant counts of $\Mm{1}{n}$ for $n \leq 10$ and hence we can compute \eqref{eq-noncubtrans} for all partitions $\lambda$ such that $\abs{\lambda} \leq 7$. 

\begin{exa} 
If we choose $\mu=[1^3]$ and $\lambda=[1^1,2^2]$ and compute formula \eqref{eq-noncubtrans}, using the results of Section \ref{sec-cubics}, we get
\begin{multline*} \frac{1}{6} \bigl(\abs{\Mmt{1}{3}^{F \cdot \sigma_{[1^3]}}} \cdot \abs{\Pbk^1([1^1,2^2])} + \abs{\Mmt{1}{4}^{F \cdot \sigma_{[1^4]} }} \cdot \abs{\Pbk^1([2^2])} + 2 \abs{\Mmt{1}{5}^{F \cdot \sigma_{[1^3,2^1]} }} \cdot \abs{\Pbk^1([1^1,2^1])}\\+2\abs{\Mmt{1}{6}^{F \cdot \sigma_{[1^4,2^1]} }} \cdot \abs{\Pbk^1([2^1])}+\abs{\Mmt{1}{7}^{F \cdot \sigma_{[1^3,2^2]} }} \cdot \abs{\Pbk^1([1^1])}+\abs{\Mmt{1}{7}^{F \cdot \sigma_{[1^4,2^2]} }} \cdot \abs{\Pbk^1(\emptyset)} \bigr) = \\ = \frac{1}{6}(8q^8-4q^7-48q^6+52q^5+65q^4-72q^3-43q^2+6q).
\end{multline*}
\end{exa}

\subsection{Quartic with three singularities} \label{sec-quart3}
Since the quartic has arithmetic genus~$3$, all of the three singularities will have $\delta =1$. Fix a type, which we will denote by $\epsilon$, a quartic of this type and a normalization morphism of the quartic among the $\PGLk{2}$ possible ones. The inverse image of the singularities will give a $\lambda$-set $S=\{P,Q,R\}$ of disjoint sets $P=\{p_1,p_2\}$, $Q=\{q_1,q_2\}$ and $R=\{r_1,r_2\}$, of one or two points on $\Pbk^1$, where $\lambda=\sum_{i=1}^3 [i^{\#c^i+\#n_1^i+\#n_2^i}]$. More precisely, $S$ contains, for $1 \leq i \leq 3$, 
\begin{itemize}  
\item $\#c^i$ $i$-sets, of sets consisting of one point; 
\item $\#n_1^i$ $i$-sets, of unordered pairs of points, that together form an $[i^2]$-set of points;
\item$\#n_2^i$ $i$-sets, of unordered pairs of points, that together form a $2i$-set of points.
\end{itemize}
The set $\{p_1,p_2,q_1,q_2,r_1,r_2\}$ will 
be a $\mu$-set where $\mu:=\sum_{i=1}^3 [i^{\#c^i+2\#n_1^i},(2i)^{\#n_2^i}].$

Let us find in how many ways we can choose $P=\{p_1,p_2\}$, $Q=\{q_1,q_2\}$ and $R=\{r_1,r_2\}$ on $\Pbk^1$ fulfilling the criteria just described. One can construct $i$ distinct $i$-sets of pairs of points out of two $i$-sets of points, and also that one can only construct one $i$-set of pairs of points out of a $2i$-set of points. From this it follows that there are $\abs{\Pbk^1(\mu)}/(v_3(\epsilon) w_3(\epsilon))$ choices of $P$, $Q$ and $R$, where 
\begin{equation}
v_j(\epsilon) := \prod_{i=1}^j i^{\#c^i} \#c^i! \cdot  i^{\#n_1^i}\#n_1^i! \cdot  i^{\#n_2^i}\#n_2^i! \;\; \text{and} \;\; w_j(\epsilon):=\prod_{i=1}^j 2^{\#n_1^i+\#n_2^i}. \end{equation} 

Fix $P$, $Q$ and $R$ as above. B\'ezout's Theorem shows that on a quartic, three singularities cannot lie on a line, and hence the following construction will give all possible normalization morphisms. Define the linear system $\mathfrak{d}$, spanned by the divisors $p_1+p_2+q_1+q_2$, $p_1+p_2+r_1+r_2$ and $q_1+q_2+r_1+r_2$. Even though the individual divisors need not be defined over $k$ they are as a triple, and hence so is $\mathfrak{d}$. Since each pair of the divisors defining $\mathfrak{d}$ has two points in common and all three have none in common, $\mathfrak{d}$ has dimension $2$ and is base point free. 

\begin{lem} \label{clm-quart3-deg1} If the morphisms induced by $\mathfrak{d}$ have degree $1$, the image will be a quartic curve with three singularities.
\end{lem}
\begin{pf} The line in $\mathfrak{d}$ passing through the divisors $p_1+p_2+q_1+q_2$ and $p_1+p_2+r_1+r_2$ has the point(s) of $P$ as base point(s). There is therefore a $\Pb^1$ of lines passing through the image of $P$, which shows that the image of $P$ is a point. 
(This argumentation also holds when the degree of the morphisms is $2$.) Since each line intersects the rest of the image quartic in at most two more points the image of $P$ must be a singularity, and the same holds for $Q$ and $R$. Since there cannot be more than three singularities on a quartic we are done. \end{pf}

Define $\mathfrak{h}$ to be the linear system spanned by the divisors $p_1+p_2$, $q_1+q_2$ and $r_1+r_2$. This is a base point free linear system defined over $k$.

\begin{lem} \label{clm-quart3-deg2} The morphisms induced by $\mathfrak{d}$ will have degree $2$ precisely if $\mathfrak{h}$ has dimension $1$. 
\end{lem}
\begin{pf} If the morphisms induced by $\mathfrak{d}$ have degree $2$ then, as noted in the proof of Lemma \ref{clm-quart3-deg1}, the image of $P$ is a point. The tangent line through the image of $P$ does not intersect the image conic in any other point and hence the corresponding divisor is $2p_1+2p_2$. The line of divisors in $\mathfrak{d}$ through $p_1+p_2+r_1+r_2$ and $p_1+p_2+q_1+q_2$ therefore contains $2p_1+2p_2$ which is equivalent to the fact that $\mathfrak{h}$ has dimension $1$. 

For the other direction, say that $\mathfrak{h}$ has dimension $1$. It follows, as was just mentioned, that $\mathfrak{d}$ contains the divisors $2p_1+2p_2$, $2q_1+2q_2$ and $2r_1+2r_2$. Assume that the image is a quartic. The line corresponding to the divisor $2p_1+2p_2$ would then intersect the image of $P$, which is a singularity, with multiplicity $4$ and hence it would have $\delta \geq 2$. Since this is also true for $Q$ and $R$ we get a contradiction. 
\end{pf}

\begin{lem} \label{clm-separable} 
If the characteristic is even and $P$, $Q$ and $R$ all consist of one point (that is, in the case of three cusps) the morphisms induced by $\mathfrak{h}$ are purely inseparable. In all other cases the morphisms induced by $\mathfrak{h}$ are separable. 
\end{lem}
\begin{pf} Any morphism between curves can be factored into a composition of separable and purely inseparable ones. The purely inseparable morphisms are (up to automorphism) compositions of the $\bar k$-linear Frobenius morphism. The $\bar k$-linear Frobenius morphism has degree equal to the characteristic and it is the identity on the topological spaces. Since the divisors of $\mathfrak{h}$ have degree $2$, we conclude that the morphisms of $\mathfrak{h}$ are separable if the characteristic is odd.

If the characteristic is even, we see that the linear system $\mathfrak{f}$ corresponding to the $\bar k$-linear Frobenius consists of the divisors of the form $2s$ for all $s$ on $\Pbk^1$. If any of $P$, $Q$ and $R$ consists of more than one point we then directly see that $\mathfrak{h}$ is distinct from $\mathfrak{f}$. On the other hand, if $P$, $Q$ and $R$ all consist of one point we see that $\mathfrak{h}$ is spanned by three distinct points of $\mathfrak{f}$. Since $\mathfrak{f}$ has dimension $1$ they must be equal. \end{pf}

\begin{lem} \label{lem-quart3-spec} In even (positive) characteristic there are no irreducible quartics with three cusps.
\end{lem}
\begin{pf} In the proof of Lemma \ref{clm-separable} we saw that in this case $\mathfrak{h}$ is equal to $\mathfrak{f}$ which has dimension $1$. We can then conclude from Lemma \ref{clm-quart3-deg2}.\end{pf}

Assume now that we are not in the case of even characteristic with three cusps. The question that remains is then for how many of the choices of $P$, $Q$ and $R$ we have that $\mathfrak{h}$ has dimension $1$. When $\mathfrak{h}$ has dimension $1$ we saw in Lemma \ref{clm-separable} that it defines $\PGLk{2}$ separable degree $2$ morphisms from $\Pbk^1$ to $\Pbk^1$ with the property that there are three points on the image curve whose preimages are $P$, $Q$ and $R$ respectively. 

\begin{dfn}[{\cite[Def. 12.1]{Jonas1}}] Let $C_{\varphi}$ be an integral non-singular projective curve over $k$ together with a separable degree $2$ morphism $\varphi$ to $\Pbk^1$. We then define
$$b_i(C_{\varphi}):=\abs{\{P \in \Pbk^1([i^1]) : \abs{\varphi^{-1}(P)}=2i, \; \varphi^{-1}(P)\subseteq C(k_i) \}},$$
$$c_i(C_{\varphi}):=\abs{ \{P \in \Pbk^1([i^1]) : \abs{\varphi^{-1}(P)}=2i, \; \varphi^{-1}(P)\nsubseteq C(k_i) \}}$$
and $r_i(C_{\varphi}):=b_i(C_{\varphi})+c_i(C_{\varphi}).$
\end{dfn}

In this terminology, what we wish to compute is the sum of 
\begin{multline*}
\phi_{\epsilon}(C_{\varphi}):=\frac{1}{v_3(\epsilon)} \cdot  \prod_{i=1}^3 \prod_{j_1=0}^{\# n_1^i-1}(b_i(C_{\varphi})-i \cdot j_1) \prod_{j_2=0}^{\# n_2^i-1}(c_i(C_{\varphi})-i \cdot j_2) \cdot \\ \cdot \prod_{j_3=0}^{\# c^i-1} (\abs{\Pbk^1([i^1])}-b_i(C_{\varphi})-c_i(C_{\varphi})-i \cdot j_3)
\end{multline*}
over representatives $\varphi$ of the linear systems of separable degree $2$ morphisms from $\Pbk^1$ to $\Pbk^1$. 

\begin{rmk} There are $q^2$ linear systems of separable degree $2$ morphisms from $\Pbk^1$ to $\Pbk^1$ when the characteristic is odd, and $q^2-1$ when the characteristic is even. This is also the number of nontrivial involutions of $\Pbk^1$ because each such gives rise to a linear system as above, and vice versa.
\end{rmk}

To do the computation we will use the definitions and results of \cite{Jonas1}, see especially the first appendix. The linear systems of separable degree $2$ morphisms from $\Pbk^1$ to $\Pbk^1$ are induced by elements of the set $P_0$, see Definition 3.1 in \cite{Jonas1} for odd characteristic and Definition 8.1 in \cite{Jonas1} for even characteristic. For an element $f$ in $P_0$ in odd characteristic (resp. $(h,f)$ in even characteristic), we write $C_f$  (resp. $C_{(h,f)}$) for the corresponding curve together with its degree $2$ morphism to $\Pb^1$. After summing $\phi_{\epsilon}$, over the curves induced by elements of $P_0$, we will not divide by $\PGLk{2}$ but by the number of elements of the group $G$ in odd characteristic and $G_0$ in even characteristic, to take into account that every linear system appears several times in $P_0$.

We have shown the following lemma.

\begin{lem} \label{lem-quart3} The number of irreducible quartics with three singularities, of a type $\epsilon$ (except for the case of three cusps in even characteristic), is equal to
\begin{equation} \label{eq-quart3} 
\frac{\abs{\Pbk^1(\mu)}}{v_3(\epsilon) w_3(\epsilon)} \cdot  
\frac{1}{\PGLk{2}}  - \begin{cases} \frac{1}{\abs{G}} \cdot \sum_{f \in P_0} \phi_{\epsilon}(C_f) & \text{if} \; \mathrm{char}(k) \neq 2; \\ \frac{1}{\abs{G_0}} \cdot \sum_{(h,f) \in P_0} \phi_{\epsilon}(C_{(h,f)}) & \text{if} \; \mathrm{char}(k) = 2. \end{cases} 
\end{equation}
\end{lem}

Equation \eqref{eq-quart3} is computable for every choice of $\epsilon$ following the arguments of Section 6 in \cite{Jonas1}. Since the application of this lemma is not immediate we present the results in Table \ref{tab-three-irred}. 

\begin{exa} 
Consider the quartic curves with one cusp over $k$ and two non-split nodes over $k$, which corresponds to the partitions $\lambda=[1^3]$ and $\mu=[1^1,2^2]$. These curves have the same number of points as $\Pbk^1$ plus two points over $k$ minus two conjugate pairs of points. Using the terminology of \cite{Jonas1} we can express the part of equation \eqref{eq-quart3} coming from the degree $2$ morphisms as
$$\frac{1}{8}\bigl(-u_0^{(1^2,1^2,1^2)}+qu_0^{(1^2,1^2)}-q u_0^{(1^2)}\bigr)$$
which also can be computed using the techniques of \cite{Jonas1}. The result of equation \eqref{eq-quart3} is then found to be
$$\frac{1}{8}(q+1)(q^2-q)(q^2-q-2)\cdot \frac{1}{\PGLk{2}}- \frac{1}{8} \begin{cases}q-3 & \text{if} \; \mathrm{char}(k) \neq 2;\\q-2 & \text{if}\; \mathrm{char}(k) = 2. \end{cases} $$
\end{exa}

\begin{table} \caption{Quartics with three singularities} \label{tab-three-irred}
\centerline{
\vbox{
\offinterlineskip
\hrule
\halign{&\vrule#& \quod #\hfil  \strut &\vrule#& \quod \hfil#\hfil \strut  &\vrule#& \quod \hfil#\hfil \strut  \cr
height2pt&\omit&&\omit&&\omit& \cr
& \bf{Description} && \bf{\# (odd char.)} && \bf{\# (even char.)} & \cr
height2pt&\omit&&\omit&&\omit&\cr
\noalign{\hrule}
height2pt&\omit&&\omit&&\omit&\cr
&\emph{\small Quartic $\mathcal Q$, sing. in $\scriptstyle [1^3]$} && && &\cr 
height2pt&\omit&&\omit&&\omit&\cr
&$\scriptstyle \star \,\,\mathcal Q=3c,2c+n_1,2c+n_2,c+2n_1$ &&$\scriptstyle \frac{1}{6},\frac{1}{4}(q-3),\frac{1}{4}(q-1),\frac{1}{8}(q-3)^2$ && $\scriptstyle 0,\frac{1}{4}(q-2),\frac{1}{4}q,\frac{1}{8}(q-2)(q-4) $ &\cr
height2pt&\omit&&\omit&&\omit&\cr
&$\scriptstyle \star \,\,\mathcal Q=c+n_1+n_2,c+2n_2$ &&$\scriptstyle \frac{1}{4}(q-1)^2,\frac{1}{8}(q-1)^2$ && $\scriptstyle \frac{1}{4}q(q-2),\frac{1}{8}q(q-2)$ &\cr
height2pt&\omit&&\omit&&\omit&\cr
&$\scriptstyle \star \,\,\mathcal Q=3n_1$ &&$\scriptstyle \frac{1}{48}(q-3)(q^2-7q+11)$ && $\scriptstyle \frac{1}{48}(q-2)(q-4)^2$ &\cr
height2pt&\omit&&\omit&&\omit&\cr
&$\scriptstyle \star \,\,\mathcal Q=2n_1+n_2$ &&$\scriptstyle \frac{1}{16}(q-1)(q^2-3q+1)$ && $\scriptstyle \frac{1}{16}q(q-2)^2$ &\cr
height2pt&\omit&&\omit&&\omit&\cr
&$\scriptstyle \star \,\,\mathcal Q=n_1+2n_2,3n_2$ &&$\scriptstyle \frac{1}{16}(q^3-2q^2-1),\frac{1}{48}(q-1)(q^2-3q+1)$ && $\scriptstyle \frac{1}{16}q^2(q-2),\frac{1}{48}q(q-2)^2$ &\cr
height2pt&\omit&&\omit&&\omit&\cr
&\emph{\small Quartic $\mathcal Q$, sing. in $\scriptstyle [1^1,2^1]$} && && &\cr 
height2pt&\omit&&\omit&&\omit&\cr
&$\scriptstyle \star \,\,\mathcal Q=c+c^2,n_1+c^2,n_2+c^2$ &&$\scriptstyle \frac{1}{2},\frac{1}{4}(q-1),\frac{1}{4}(q-3)$ && $\scriptstyle 0,\frac{1}{4}q,\frac{1}{4}(q-2)$ &\cr
height2pt&\omit&&\omit&&\omit&\cr
&$\scriptstyle \star \,\,\mathcal Q=c+n_1^2,c+n_2^2$ &&$\scriptstyle \frac{1}{4}(q^2-2q-1),\frac{1}{4}(q+1)(q-1) $ && $\scriptstyle \frac{1}{4}q(q-2),\frac{1}{4}q^2$ &\cr
height2pt&\omit&&\omit&&\omit&\cr
&$\scriptstyle \star \,\,\mathcal Q=n_1+n_1^2,n_2+n_1^2$ &&$\scriptstyle \frac{1}{8}(q-1)(q^2-q-1),\frac{1}{8}(q-3)(q^2-q-3)$ && $\scriptstyle \frac{1}{8}q^2(q-2),\frac{1}{8}(q-2)(q^2-2q-4)$ &\cr
height2pt&\omit&&\omit&&\omit&\cr
&$\scriptstyle \star \,\,\mathcal Q=n_1+n_2^2,n_2+n_2^2$ &&$\scriptstyle \frac{1}{8}(q+1)(q^2-q+1),\frac{1}{8}(q-1)(q^2-q-1)$ && $\scriptstyle \frac{1}{8}q^3,\frac{1}{8}q^2(q-2)$ &\cr
height2pt&\omit&&\omit&&\omit&\cr
&\emph{\small Quartic $\mathcal Q$, sing. in $\scriptstyle [3^1]$} && && &\cr 
height2pt&\omit&&\omit&&\omit&\cr
&$\scriptstyle \star \,\,\mathcal Q=c^3,n_1^3$ &&$\scriptstyle \frac{1}{3},\frac{1}{6}(q^3-q^2-q-3)$ && $\scriptstyle 0,\frac{1}{6}(q-2)(q^2+q+1)$ &\cr
height2pt&\omit&&\omit&&\omit&\cr
&$\scriptstyle \star \,\,\mathcal Q=n_2^3$ &&$\scriptstyle \frac{1}{6}(q-1)(q^2+1)$ && $\scriptstyle \frac{1}{6}q(q^2-q+1)$ &\cr
height2pt&\omit&&\omit&&\omit&\cr
} \hrule}}
\end{table}

\section{Two singularities} \label{sec-two}
The types of degree $4$ curves with two singularities and irreducible components of degree at most $2$, are ones consisting of two conics or one conic with two lines.

\begin{table}[htbp] \caption{Two singularities} \label{tab-two-red}
\centerline{
\vbox{
\offinterlineskip
\hrule
\halign{&\vrule#& \quod #\hfil \quod \strut &\vrule#& \quod \hfil#\hfil \strut \quod \cr
height2pt&\omit&&\omit& \cr
& \bf{Description} && \bf{\#} & \cr
height2pt&\omit&&\omit&\cr
\noalign{\hrule}
height2pt&\omit&&\omit&\cr
&\emph{\small Conic + tgt line + trans. line, int. in $\scriptstyle [1^2]$} &&$\scriptstyle \frac{1}{q-1}$ &\cr
height2pt&\omit&&\omit&\cr
&\emph{\small $\lambda$-set of conics, int. in $\mu$-set} && &\cr
height2pt&\omit&&\omit&\cr
&$\scriptstyle \star \,\,(\lambda,\mu)=([1^2],[1^2]),([1^2],[2^1]),([2^1],[1^2]),([2^1],[2^1])$ &&$\scriptstyle \frac{3q-4}{4(q-1)},\frac{q-2}{4(q+1)},\frac{3q-2}{4(q-1)},\frac{q}{4(q+1)}$&\cr
height2pt&\omit&&\omit&\cr
} \hrule}}
\end{table}

\subsection{Singular cubic with a line intersecting in one smooth point on the cubic} \label{sec-cubone}
In the same way as in Sections \ref{sec-cubtrans} and \ref{sec-cubtwo}, but where $Q=(q_1)$ is a $[1^1]$-tuple and $q_1=q_2=q_3$, we show the following lemma.

\begin{lem} \label{lem-cubone} The number of singular cubics with a line intersecting in one smooth point on the cubic, of a type determined by the partition $\lambda$, is equal to 
$$ \frac{1}{\zee{\lambda}} \abs{\Pbk^1(\lambda)} \cdot \abs{X([1^1])} \cdot \frac{1}{\PGLk{2}} \;\; \text{where} \;\; X=\Pbk^1 \setminus \{a \; \lambda-set \}.$$
\end{lem}

\subsection{Singular cubic with a line through the singularity intersecting in one more point} \label{sec-cubsingone}
The inverse image of the singularity under a normalization morphism of the singular cubic gives a $\lambda$-set $P$ where $\abs{\lambda}=1$ or $2$. 

\begin{lem} \label{lem-cubsingone} The number of singular cubics with a line through the singularity intersecting in one more point, of a type determined by the partition $\lambda$, is equal to 
$$\frac{1}{\zee{\lambda}}\abs{\Pbk^1(\lambda)} \cdot (\abs{X([1^2])}+\abs{X([1^1])}) \cdot \frac{1}{\PGLk{2}} \;\; \text{where} \;\; X=\Pbk^1 \setminus \{a \; \lambda-set \}. $$
\end{lem}
\begin{pf} Let $q$ be the smooth point of intersection with the line. The tangent line through $q$ intersects the singular cubic in either one more smooth point or no other point. The singular cubic together with this tangent line is therefore a curve of a type found in either Section \ref{sec-cubtwo} or \ref{sec-cubone}, and this construction gives a bijection. The result now follows from Lemmas \ref{lem-cubtwo} and \ref{lem-cubone}. 
\end{pf}

\subsection{Non-singular cubic with a line intersecting in two points} \label{sec-noncubtwo}
In this case we can mimic Section \ref{sec-noncubtrans} with the difference that $P=(p_1,p_3)$ is a $[1^2]$-tuple and $p_1=p_2$.

\begin{lem} \label{lem-noncubtwo} The number of non-singular cubics with a line intersecting in two smooth points on the cubic, of a type determined by the partition $\mu$, and with a $\lambda$-tuple of points is equal to
\begin{equation}  \label{eq-noncubtwo} 
\sum_{i=0}^{\abs{\lambda}} \sum_{\substack{\kappa \vdash i \\ \kappa \leq \lambda}}\frac{\zee{\lambda}}{\zee{\kappa} \zee{\lambda-\kappa}} \abs{\Mmt{1}{2+i}^{F \cdot \sigma_{\mu+\kappa}} } \cdot \abs{\Pbk^1(\lambda-\kappa)}.  \end{equation}
\end{lem}

The results of Section \ref{sec-cubics} will enable us to compute equation \eqref{eq-noncubtwo} for all partitions $\lambda$ where $\abs{\lambda} \leq 8$.

\subsection{Quartic with one $\delta = 2$ singularity and one $\delta =1$ singularity} \label{sec-quart12} 
Let us recall Definition \ref{dfn-sings}. The inverse image of the singularity with $\delta =2$ (resp. $\delta =1$), under any normalization morphism, will be a $\lambda$-set $P$ (resp. $\mu$-set $Q$) on the normalization $\Pbk^1$, where $\abs{\lambda}$ (resp. $\abs{\mu}$) is $1$ or $2$ . The sets $P$ and $Q$ are disjoint.

Fix $P=\{p_1,p_2\}$ and $Q=\{q_1,q_2\}$ fulfilling the conditions above and consider the $4$-dimensional space of divisors $|2p_1+2p_2|$. It contains the point $2p_1+2p_2$ that lies outside the plane of divisors of the form $q_1+q_2+r_1+r_2$ for any points $r_1$ and $r_2$ on $\Pbk^1$. To get a linear system of dimension $2$, we wish to choose a line in this plane that passes through the divisor $q_1+q_2+p_1+p_2$. We should avoid the two (possibly equal) lines consisting of divisors of the form $q_1+q_2+p_1+s_1$ for any $s_1$ on $\Pbk^1$ and $q_1+q_2+p_2+s_2$ for any $s_2$ on $\Pbk^1$. This makes the linear system base point free, because the divisors on the chosen line do not all contain either $p_1$ or $p_2$. Moreover, we should also avoid the line through the divisor $2q_1+2q_2$, which is distinct from the previous lines. Thus, we have $q- \lambda_1 $ possibilities to choose $\mathfrak{d}$. 

\begin{lem} \label{clm-quart12} The image of any morphism induced by $\mathfrak{d}$ is a quartic with one $\delta =2$ singularity and one $\delta =1$ singularity.
\end{lem}
\begin{pf} By construction, there is a $\Pb^1$ of lines through the image of $Q$ and hence it must consist of one point. If the image curve is a conic, the divisor corresponding to the tangent line through the image of $Q$ is equal to $2q_1+2q_2$. But this divisor is not in $\mathfrak{d}$ and hence the image must be a quartic and $Q$ a singularity. The two lines corresponding to the divisors $2p_1+2p_2$ and $q_1+q_2+p_1+p_2$ both pass through the image of $P$ and hence it must also be a singularity. Furthermore, the line corresponding to $2p_1+2p_2$ intersects the singularity with multiplicity $4$ and hence the image of $P$ must have $\delta =2$ and the image of $Q$ must have $\delta =1 $. 
\end{pf} 

\begin{lem} \label{lem-quart12} The number of quartics with one $\delta =2$ singularity and one $\delta =1$ singularity, of a type determined by the partitions $\lambda$ and $\mu$, is equal to 
$$\frac{1}{\zee{\lambda}\zee{\mu}} \abs{\Pbk^1(\lambda)} \cdot \abs{X(\mu)} \cdot (q-\lambda_1) \cdot \frac{1}{\PGLk{2}} \;\; \text{where} \;\; X=\Pbk^1 \setminus \{a \; \lambda-set \}. $$
\end{lem}

\subsection{Quartic with two $\delta =1$ singularities} \label{sec-quart2}
Fix a type of quartic with two $\delta =1$ singularities, which we denote by $\epsilon$, and a curve over $k$ of this type. The normalization $C$ of this curve has genus $1$. The inverse images under a normalization morphism of each of the singularities give a $\lambda$-set $S=\{P,Q\}$ of disjoint sets $P=\{p_1,p_2\}$ and $Q=\{q_1,q_2\}$ of one or two points on $C$, where $\lambda:=\sum_{i=1}^2 [i^{\#c^i+\#n_1^i+\#n_2^i}]$. The set $\{p_1,p_2,q_1,q_2\}$ will be a $\mu$-set where $\mu:=\sum_{i=1}^2 [i^{\#c^i+2\#n_1^i},(2i)^{\#n_2^i}].$ Moreover, if $\lambda=[1^2]$, the partition $\mu$ will be a member of the set $\Delta_1:=\{[1^2],[1^3],[1^1,2^1],[1^4],[1^2,2^1],[2^2] \}$ and if $\lambda=[2^1]$ it will be a member of $\Delta_2:=\{[2^1],[2^2],[4^1]\}.$ In the same way as in Section \ref{sec-quart3} we then find that there are $\abs{C(\mu)}/(v_2(\epsilon) w_2(\epsilon))$ choices of $P$ and $Q$. 

\begin{lem} \label{clm-quart2-type} The type $\epsilon$ is determined by $\lambda$ and $\mu$.
\end{lem}
\begin{pf} If $\lambda=[1^2]$ we have $\# c + \#n_1 +\#n_2=\lambda_1$, $\# c + 2\cdot \#n_1=\mu_1$ and $\#n_2=\mu_{2}$, and if $\lambda=[2^1]$ we have $\# c^2 + \#n_1^2 +\#n_2^2=\lambda_2$, $\# c^2 + 2\cdot \#n_1^2=\mu_2$ and $\#n_2^2=\mu_{4}$. \end{pf}

Fix a non-singular curve $C$ of genus $1$ over $k$ and choose points $P=\{p_1,p_2\}$ and $Q=\{q_1,q_2\}$ on $C$ as above. It follows from Riemann-Roch that the complete linear system $|p_1+p_2+q_1+q_2|$ has dimension $3$. This theorem also shows that the divisors in $|p_1+p_2+q_1+q_2|$ of the form $p_1+p_2+s_1+s_2$, for any $s_1$ and $s_2$ on $C$, form a line. Similarly, the divisors in the linear system of the form $q_1+q_2+s_1+s_2$, for any $s_1$ and $s_2$ on $C$, form a line which is distinct from the previous one. Define the linear system $\mathfrak{d}$ to be the the plane inside $|p_1+p_2+q_1+q_2|$ that contains these two lines. It is base point free.

\begin{lem} \label{clm-quart2-deg2} The morphisms induced by $\mathfrak{d}$ will have degree $2$ precisely when the divisors $p_1+p_2$ and $q_1+q_2$ are linearly equivalent.
\end{lem}
\begin{pf} Assume first that the morphisms induced by $\mathfrak{d}$ have degree $2$. Since there is a $\Pb^1$ of divisors which have the point(s) of $P$ as base point(s), the image of $P$ is a point. The tangent line through the image of $P$ corresponds to the divisor $2p_1+2p_2$ which then lies in $|p_1+p_2+q_1+q_2|$. This implies that the divisors $p_1+p_2$  and $q_1+q_2$ are linearly equivalent.

Suppose that $p_1+p_2$ and $q_1+q_2$ are linearly equivalent. The divisors $2p_1+2p_2$ and $2q_1+2q_2$ are then found in $\mathfrak{d}$. The lines corresponding to these divisors intersect the images of $P$ and $Q$ respectively in no other points. If the image curve would be a quartic then $P$ and $Q$ would be singularities with $\delta \geq 2$, which is not possible. 
\end{pf}

\begin{lem} \label{clm-quart2-deg1} If the morphisms induced by $\mathfrak{d}$ have degree $1$, the image will be a quartic curve with two $\delta =1$ singularities.
\end{lem}
\begin{pf} The existence of a $\Pb^1$ of divisors with the point(s) of $P$ (resp. $Q$) as base point(s) implies that the image of $P$ (resp. $Q$) will be a singularity. If $P$ would have $\delta =2$, then the line in the tangent cone would correspond to the divisor $2p_1+2p_2$ which would lie in $\mathfrak{d}$. But this is not possible by Lemma \ref{clm-quart2-deg2}. 
\end{pf}

If $p_1+p_2$ and $q_1+q_2$ are linearly equivalent the complete linear system $|p_1+p_2|=|q_1+q_2|$, which we denote by $\mathfrak{h}$, has dimension $1$ and is, by Riemann-Roch, base point free. 

\begin{lem} \label{clm-quart2-separable} Suppose that $p_1+p_2$ and $q_1+q_2$ are linearly equivalent. Then the  morphisms induced by $\mathfrak{h}$ are separable. 
\end{lem}
\begin{pf}
Since the morphisms of $\mathfrak{h}$ have degree $2$ they must be either separable or purely inseparable. But the genera of the source and image of a purely inseparable morphism are the same. 
\end{pf}

Lemmas \ref{clm-quart2-deg2} and \ref{clm-quart2-separable} tell us that $p_1+p_2$ and $q_1+q_2$ are linearly equivalent precisely if there is a separable morphism of degree $2$ from $C$ to $\Pbk^1$ and a $\lambda$-set of points on the image $\Pbk^1$ whose preimages are equal to $P$ and $Q$ respectively. The number of choices of $P$ and $Q$, for which $p_1+p_2$ and $q_1+q_2$ are linearly equivalent, is therefore equal to the sum of 
\begin{multline*}
\phi_{\lambda,\mu}(C_{\varphi}):= \frac{1}{v_2(\epsilon)} \cdot \prod_{i=1}^2 \Bigl( \prod_{j_1=0}^{\# n_1^i-1}(b_i(C_{\varphi})-i \cdot j_1) \prod_{j_2=0}^{\# n_2^i-1}(c_i(C_{\varphi})-i \cdot j_2) \cdot \\ \cdot \prod_{j_3=0}^{\# c^i-1} (\abs{\Pbk^1([i^1])}-b_i(C_{\varphi})-c_i(C_{\varphi})-i \cdot j_3) \Bigr)
\end{multline*}
over representatives $\varphi$ of the linear systems of separable degree $2$ morphisms from $C$ to $\Pbk^1$. 

Fix a partition $\kappa=[1^{\kappa_1},\ldots,\nu^{\kappa_{\nu}}]$. We are ready to find a formula for the number of quartics with two singularities, of a type $\epsilon$ determined by $\lambda$ and $\mu$, together with a $\kappa$-tuple of points. We should take into account that the $\kappa$-tuple may contain one or both of the singularities. There are two parts to this computation. By Lemma~\ref{clm-quart2-deg1} and by the arguments of Section \ref{sec-quartics}, we should first count the number of $k$-isomorphism classes of non-singular curves $C$ of genus $1$ over $k$ together with a $\mu$-set of points on $C$, and a $\kappa$-tuple of points on the union of $C \setminus \{a \; \mu-set \}$ and a $\lambda$-set. As usual, the pointed curves should be counted with weight the reciprocal of their number of $k$-automorphisms. This count can be formulated in terms of $\s_n$-equivariant counts of $\Mm{1}{n}$ where $n \leq \abs{\mu}+\abs{\kappa}$. We should then subtract the number of cases when the linear system $\mathfrak{d}$ induces morphisms of degree $2$. By the arguments directly above, this number is equal to the sum of
\begin{equation*}\psi^{\kappa}_{\lambda,\mu}(C_{\varphi}):=\phi_{\lambda,\mu}(C_{\varphi}) \cdot \abs{X_{\lambda,\mu}^{\varphi}(\kappa)}  \;\; \text{where} \;\; X_{\lambda,\mu}^{\varphi}=\{a \; \lambda-set \} \amalg C_{\varphi} \setminus \{a \; \mu-set \},
\end{equation*}
over representatives $\varphi$ of the $k$-isomorphism classes of linear systems of separable degree $2$ morphisms from a genus $1$ curve $C$ to $\Pbk^1$. This should also be counted with weight the reciprocal of their number of $k$-automorphisms. As in Section \ref{sec-quart3} we will use the results of \cite{Jonas1}, and realize the linear systems of separable degree $2$ morphisms from $C$ to $\Pb^1_k$ as elements of the set $P_1$. 

Actually, we are not able to perform the second part of this computation for partitions $\kappa$ with $\abs{\kappa} > 3$ using the results of \cite{Jonas1}. Nor have we made $\s_{11}$-equivariant counts of $\Mm{1}{11}$ which are necessary for the first part of the computation when $\abs{\mu}=4$ and $\abs{\kappa}=7$. But we will see that these problems dissolve when for each $\lambda$ we sum the contributions over $\mu$. 

\begin{lem} \label{lem-quart2-k} The number of irreducible quartics with two $\delta =1$ singularities defined over $k$, of any type $\epsilon$, together with a $\kappa$-tuple of points is equal to
\begin{multline} \label{eq-quart2-k} 
\sum_{i=0}^2 \frac{2!}{(2-i)!} \binom{\kappa_1}{i} \sum_{\mu \in \Delta_1} \frac{1}{v_2(\epsilon) w_2(\epsilon)} \cdot \abs{\Mmt{1}{\abs{\mu}+\abs{\kappa}-i}^{F \cdot \sigma_{\mu+\kappa-[1^i]}}}\\
-\begin{cases} \frac{1}{\abs{G}}  \sum_{f \in P_1} \sum_{\mu \in \Delta_1} \psi^{\kappa}_{[1^2],\mu}(C_f) & \text{if} \; \mathrm{char}(k) \neq 2; \\ \frac{1}{\abs{G_1}} \sum_{(h,f) \in P_1} \sum_{\mu \in \Delta_1} \psi^{\kappa}_{[1^2],\mu}(C_{(h,f)}) &  \text{if} \; \mathrm{char}(k) = 2. \end{cases} 
\end{multline}
\end{lem}

\begin{lem} \label{lem-quart2-k2} The number of irreducible quartics with a conjugate pair of $\delta =1$ singularities, of any type $\epsilon$, together with a $\kappa$-tuple of points is equal to
\begin{multline} \label{eq-quart2-k2} 
\sum_{i=0}^1 2^i\kappa_2^i \sum_{\mu \in \Delta_2}\frac{1}{v_2(\epsilon) w_2(\epsilon)} \cdot \abs{\Mmt{1}{\abs{\mu}+\abs{\kappa}-2i}^{F \cdot \sigma_{\mu+\kappa-[2^i]}}} \\
-\begin{cases} \frac{1}{\abs{G}}  \sum_{f \in P_1} \sum_{\mu \in \Delta_2} \psi^{\kappa}_{[2^1],\mu}(C_f) & \text{if} \; \mathrm{char}(k) \neq 2; \\ \frac{1}{\abs{G_1}} \sum_{(h,f) \in P_1} \sum_{\mu \in \Delta_2} \psi^{\kappa}_{[2^1],\mu}(C_{(h,f)}) &  \text{if} \; \mathrm{char}(k) = 2. \end{cases} 
\end{multline}
\end{lem}

\begin{exa} Let us count the irreducible quartics with a conjugate pair of $\delta =1$ singularities of any type, and a conjugate pair of points. With the information from Section \ref{sec-cubics} we can compute the first part of equation \eqref{eq-quart2-k2}:
\begin{multline*}
\frac{1}{2}(\abs{\Mmt{1}{4}^{F \cdot \sigma_{[2^2]}}}+2\abs{\Mmt{1}{2}^{F \cdot \sigma_{[2^1]}}})+
\frac{1}{4}(\abs{\Mmt{1}{6}^{F \cdot \sigma_{[2^3]}}}+2\abs{\Mmt{1}{4}^{F \cdot \sigma_{[2^2]}}})+\\
+\frac{1}{4}(\abs{\Mmt{1}{6}^{F \cdot \sigma_{[2^1,4^1]}}}+2\abs{\Mmt{1}{4}^{F \cdot \sigma_{[4^1]}}}) = \frac{1}{2}(q^6+q^4-2q^2-q+1).
\end{multline*}
Using the terminology and techniques of \cite{Jonas1} we can express and compute the second part of the equation:
$$\frac{1}{2}(q(q^2-q)^2+(q^2-q-2) u_1^{(2^1)}+(q^2-q) u_1^{(1^2)})= \frac{1}{2}(q^5-q^4+q^2-q-2),$$
which is independent of the characteristic.
\end{exa}

We will now 
show that equations \eqref{eq-quart2-k} and \eqref{eq-quart2-k2} are computable for all $\kappa$ with $ \abs{\kappa}\leq7$. 

\subsubsection{The first part of equations \eqref{eq-quart2-k} and \eqref{eq-quart2-k2}} \label{sec-firstpart} By the arguments of Section \ref{sec-quartics} we can show that, for any $n \geq 1$ and $\sigma \in S_n$, 
\begin{equation} \label{eq-m1n-count}
\abs{\Mmt{1}{n}^{F \cdot \sigma}} = \sum_{[(C,p)]\in \Mm{1}{1}(k)/\cong_k} \frac{\abs{C(c(\sigma))}}{\abs{\Aut_k \bigl((C,p)\bigr)} \cdot \abs{C(k)}}.
\end{equation}
Using the relation $\abs{C(k_i)}=q^i+1-a_i(C)$ and equation \eqref{eq-Xlambda} we can consider $\abs{C(c(\sigma))}$ as a polynomial in the variables $q$, of degree $0$, and $a_i(C)$, of degree $i$, for $1 \leq i \leq n$. We then find that if we have made $\s_i$-equivariant counts of $\Mm{1}{i}$ for every $1 \leq i \leq n$ we can compute
$$\sum_{[(C,p)]\in \Mm{1}{1}(k)/\cong_k} \frac{p(C)}{\abs{\Aut_k \bigl((C,p)\bigr)}}$$
for any polynomial $p(C)$ in the variables $q$ and $a_i(C)$ of degree at most $n-1$.

In our case $c(\sigma)=\mu+\kappa$, and we will consider 
$$\xi_1(C):= \sum_{\mu \in \Delta_1} \frac{\abs{C(\mu+\kappa)}}{v_2(\epsilon) w_2(\epsilon)} \;\; \text{and} \;\; \xi_2(C):= \sum_{\mu \in \Delta_2} \frac{\abs{C(\mu+\kappa)}}{v_2(\epsilon) w_2(\epsilon)}$$ 
as polynomials in $q$ and $a_i(C)$. The degree $\abs{\mu+\kappa}$ part of $\xi_1(C)$ is equal to 
$\bigl(\frac{1}{8} a_1(C)^4-\frac{1}{4} a_2(C)a_1(C)^2+\frac{1}{8}a_2(C)^2\bigr) \cdot \prod_{i=1}^{\nu} a_i^{\kappa_i}(C)$
and that of $\xi_2(C)$ is equal to 
$\bigl (\frac{1}{4} a_2(C)^2-\frac{1}{4} a_4(C)\bigr) \cdot \prod_{i=1}^{\nu} a_i^{\kappa_i}(C).$
In Section \ref{sec-cubics} we will see that for genus $1$ curves, $a_i(C)$ is expressible in terms of $q$ and $a_1(C)$. In particular, we have $a_2(C)=a_1(C)^2-2q$ and $a_4(C)=a_1(C)^4 -4qa_1(C)^2+2q^2$. This shows that, formulated in $q$ and $a_1(C)$, the degree $\abs{\mu+\kappa}$ part of $\xi_1(C)$ and $\xi_2(C)$ is equal to $0$. 

Furthermore, for any nonempty partition $\chi$, the polynomial $\abs{C(\chi)}$ formulated in $q$ and $a_1(C)$ is divisible by $\abs{C(k)}$. It follows that $\xi_1(C)/\abs{C(k)}$ and $\xi_2(C)/\abs{C(k)}$ are polynomials of degree at most $2+\abs{\kappa}$. Hence, if $\abs{\kappa} \leq 7$, the $\s_n$-equivariant counts of $\Mm{1}{n}$ for $n \leq 10$ found in Section \ref{sec-cubics}, suffice to compute the first part of equations \eqref{eq-quart2-k} and \eqref{eq-quart2-k2}.

\subsubsection{The second part of equations \eqref{eq-quart2-k} and \eqref{eq-quart2-k2}}
Consider $\psi_{\lambda,\mu}^{\kappa}(C_f)$ as a polynomial in the variables $q$, of degree $0$, and $a_i(C_f)$, $b_i(C_f)$ and $c_i(C_f)$ of degree $i$. It is equal to the degree $2$ polynomial $\phi_{\lambda,\mu}(C_f)$ in $q$, $b_i(C_f)$ and $c_i(C_f)$, times a polynomial in $q$ and $a_i(C_f)$ of degree equal to the weight of $\kappa$. 

We saw in the first appendix of \cite{Jonas1} that, after summing over all $f$ in $P_1$, we can compute all polynomials of degree at most $5$ in all the four types of variables. The results of Section \ref{sec-cubics} together with Section 3.1 in \cite{Jonas1} show that, after summing over all $f$ in $P_1$,  we can compute all polynomials of degree at most $9$ in the variables $q$ and $a_i(C_f)$. Finally, Remark 12.9 in \cite{Jonas1} tells us that the sum over all $f$ in $P_1$ of $r_1(C_f)$ times a polynomial in $q$ and $a_i(C_f)$ of odd degree is equal to $0$. 
The following lemma is therefore sufficient to conclude that we can compute the second part of equations \eqref{eq-quart2-k} and \eqref{eq-quart2-k2}. 

\begin{lem} Fix a partition $\kappa$ of weight $\leq 7$ and put 
$$\psi_1:=\sum_{\mu \in \Delta_1}\psi^{\kappa}_{[1^2],\mu}(C_f) \;\; \text{and} \;\; \psi_2:=\sum_{\mu \in \Delta_2}\psi^{\kappa}_{[2^1],\mu}(C_f).$$ 
The part of $\psi_1$ and $\psi_2$ that has degree at least $6$ and that has nonzero degree in $b_i(C_f)$ or $c_i(C_f)$ is of the form $r_1(C_f)$ times a polynomial in $q$ and $a_i(C_f)$ of degree~$5$. 
\end{lem}
\begin{pf} Let us write $\lambda-\mu=[1^{\alpha_{\lambda,\mu}},2^{\beta_{\lambda,\mu}},4^{\delta_{\lambda,\mu}}]$, where we in this case also allow negative exponents. We then have
\begin{multline*}
\psi_1=\sum_{\mu \in \Delta_1} \phi_{[1^2],\mu}(C_f) \cdot \prod_{i=1}^{\kappa_1}(q+1+\alpha_{[1^2],\mu}-a_1(C_f))\cdot \\\cdot \prod_{i=1}^{\kappa_2}(q^2-q+2\beta_{[1^2],\mu}-a_2(C_f)+a_1(C_f)) \cdot \abs{X_{\lambda,\mu}^f([3^{\kappa_3},4^{\kappa_4},\ldots]) }
\end{multline*}
and a similar expression for $\psi_2$. The equalities
$\sum_{\mu \in \Delta_1} \phi_{[1^2],\mu}(C_f)=\frac{1}{2}(q^2+q)$ and $\sum_{\mu \in \Delta_2} \phi_{[2^1],\mu}(C_f) =\frac{1}{2}(q^2-q)$ show that we need only care about the parts of $\psi_1$ and $\psi_2$ containing $\alpha_{\lambda,\mu}$ or $\beta_{\lambda,\mu}$ or $\delta_{\lambda,\mu}$. Note that $\alpha_{[2^1],\mu}=0$ for all $\mu$ in $\Delta_2$. The first expression to check is therefore a polynomial in $q$ and $a_i(C_f)$ of degree $4$, $5$ or $6$, times $\sum_{\mu \in \Delta_1} \alpha_{[1^2],\mu} \cdot \phi_{[1^2],\mu}(C_f), $ which is equal to harmless $q a_1(C_f)$. Next in line are polynomials in $q$ and $a_i(C_f)$ of degree $4$ or $5$, times $\sum_{\mu \in \Delta_1} \alpha_{[1^2],\mu}^2 \cdot \phi_{[1^2],\mu}(C_f) = a_1(C_f)^2+(q-1)r_1(C_f)$ or $\sum_{\mu \in \Delta_1} \beta_{[1^2],\mu} \cdot \phi_{[1^2],\mu}(C_f)=-\frac{1}{2}q\bigl(a_1(C_f)+r_1(C_f)\bigr)$ or $\sum_{\mu \in \Delta_2} \beta_{[2^1],\mu} \cdot \phi_{[2^1],\mu}(C_f) = \frac{1}{2} \bigl(a_2(C_f)+r_1(C_f)\bigr).$ Finally there are the ones coming from a polynomial in $q$ and $a_i(C_f)$ of degree $4$, times
$\sum_{\mu \in \Delta_1} \alpha_{[1^2],\mu}^3 \cdot \phi_{[1^2],\mu}(C_f) = a_1(C_f)(3r_1(C_f)+q-3)$ or $\sum_{\mu \in \Delta_1} \alpha_{[1^2],\mu}\cdot \beta_{[1^2],\mu} \cdot \phi_{[1^2],\mu}(C_f) =-\frac{1}{2}\bigl(a_1(C_f)^2+a_1(C_f)r_1(C_f)+(q-1)(a_1(C_f)+r_1(C_f)\bigr). $

Since this covers all the possibilities of getting a monomial of degree at least $6$ that contains $b_i(C_f)$ or $c_i(C_f)$ we are done. \end{pf}

\section{Linear subspaces} \label{sec-linsub}
In this section we will find the dimensions of all linear spaces $L_{P,S}$ for $P$ a $\lambda$-tuple where $\abs{\lambda} \leq 7$ and $S$ a $\mu$-set, where $\abs{\mu} \leq 1$. That $S$ should consist of at most one point is the consequence of the choice $M=1$ in Section \ref{sec-M1}.

\begin{dfn}
For any sets (or tuples) $Q$ and $T$ of points in $\Pb^2(\bar k)$, let $L'_{Q,T}$ be the subspace of the $\Pb^{14}(\bar k)$ of plane quartic curves over $\bar k$ that contain $Q$ and have singularities at the points of $T$.
\end{dfn}

\begin{lem} \label{lem-equaldim} 
Let $Q$ and $T$ be sets of points over $\bar k$ such that there is a line $l$ that contains $q$ points of $Q$ and $t$ points of $T$, which are all distinct. Suppose that $q+2t > 4$. Then, if $p$ is any point lying on the line $l$, we have 
$\dim L'_{Q \cup p,T} = \dim L'_{Q,T}.$ 
\end{lem}

\begin{lem} \label{lem-vectorspaces}
If $Q$ and $T$ are sets fixed by Frobenius (in the sense of Remark \ref{rmk-subschemes}) then $\mathrm{dim}_k \, L_{Q,T} = \mathrm{dim}_{\bar k} \, L'_{Q,T}. $
\end{lem}

The dimension of $L_{Q,T}$ will be called the \emph{expected dimension} if $\dim L_{Q,T}= 14-\abs{Q}-3\abs{T}$. All cases of $Q$ and $T$ for which the dimension of $L_{Q,T}$ is not the expected will be called \emph{special}. 

\begin{dfn} \label{dfn-Aa} Let the set $A_{\lambda,\mu}$ consist of all $\lambda$-tuples of points $P$ such that a $\lambda-\mu$-subtuple of points of $P$ lies on a line and the remaining $\mu$-tuple lies outside this line. If $\mu \leq \lambda$ let us also put
$$a_{\lambda,\mu}:= (q^2+q+1) \frac{\zee{\lambda}}{\zee{\mu}\zee{\lambda-\mu}}  \abs{\Pbk^1(\lambda - \mu)} \abs{X(\mu)} \;\;\text{where} \;\; X:=\Pbk^2\setminus\Pbk^1,$$
and if $\mu \nleq \lambda $ put $a_{\lambda,\mu}:=0$.
\end{dfn}

\begin{lem} \label{lem-Aa} If $\abs{\lambda-\mu} > \abs{\mu}+1$ then $\abs{A_{\lambda,\mu}}=a_{\lambda,\mu}$.\end{lem}
\begin{pf} The number $a_{\lambda,\mu}$ counts the choices of a line, a $\lambda-\mu$-set of points on the line, a $\mu$-set of points outside the line, and an ordering of the chosen $\lambda$-set of points. For a fixed $\lambda$-tuple $P$ there can, by the assumption on $\lambda$ and $\mu$, only be one $\lambda-\mu$-subtuple of $P$ that lies on a line, which will be called the \emph{distinguished line}. Hence $a_{\lambda,\mu}$ counts the members of $A_{\lambda,\mu}$. \end{pf} 

\begin{ass}
From now on $S$ will consist of one point over $k$.
\end{ass}

For simplicity we put $P_i:=\abs{\Pb^i(k)}$. 

\subsection{Weight up to $3$} \label{sec-w3} 
If $P$ is a $\lambda$-tuple and $\abs{\lambda} \leq 3$, the dimension of $L_{P,S}$ is the expected as long as $P$ and $S$ are disjoint. Since the singularity is defined over $k$ it can only be equal to one of the $\lambda_1$ points of $P$ that are defined over $k$, and hence 
$$\linsub{1}=\Pbpl \bigl(P_{14-\abs{\lambda}}- \lambda_1 P_{12-\abs{\lambda}}- (q^2+q+1-\lambda_1)P_{11-\abs{\lambda}}\bigr).$$

\subsection{Weight $4$} \label{sec-w4} 
For each of the $a_{\lambda,\emptyset}$ choices of a $\lambda$-tuple $P$ in $A_{\lambda,\emptyset}$ it follows from Lemma \ref{lem-equaldim} that $\dim L_{P,S} = 8$ if $S$ lies on the distinguished line of $P$. The only other special case is when $P$ and $S$ are not disjoint, and therefore
$$\linsub{1}=a_{\lambda,\emptyset}\bigl(P_{10}- (q+1)P_{8}- q^2 P_{7}\bigr)+\bigl(\Pbpl-a_{\lambda,\emptyset}\bigr)\bigl(P_{10}- \lambda_1 P_{8}- (q^2+q+1-\lambda_1) P_{7}\bigr). $$

\subsection{Weight $5$} \label{sec-w5} 
The only new special case is when $P$ is in $A_{\lambda,[1^1]}$ and $S$ lies on the distinguished line in $P$, and hence 
\begin{multline*} 
\linsub{1}=a_{\lambda,\emptyset}(P_{9}- (q+1)P_{8}- q^2 P_{6}) +a_{\lambda,[1^1]}(P_{9}- (q+1)P_{7}-P_{7}- (q^2-1) P_{6})\\+(\Pbpl-a_{\lambda,\emptyset}-a_{\lambda,[1^1]})(P_{9}- \lambda_1 P_{7}- (q^2+q+1-\lambda_1) P_{6}). \end{multline*}

\subsection{Weight $6$} \label{sec-w6}
If $P$ lies in $A_{\lambda,\emptyset}$ then $\dim L_{P} = 9$. The other new special cases are when $P$ lies in either $A_{\lambda,[1^2]}$ or $A_{\lambda,[2^1]}$ and $S$ lies on the distinguished line of $P$. We then conclude that 
\begin{multline*} 
\linsub{1}=a_{\lambda,\emptyset}(P_{9}- (q+1)P_{8}- q^2 P_{6}) +a_{\lambda,[1^1]}(P_{8}- (q+1)P_{7}-P_{6}- (q^2-1) P_{5})\\+a_{\lambda,[1^2]}(P_{8}- (q+1)P_{6}-2P_{6}- (q^2-2) P_{5})+a_{\lambda,[2^1]}(P_{8}- (q+1)P_{6}- q^2 P_{5}) \\+(\Pbpl-a_{\lambda,\emptyset}-a_{\lambda,[1^1]}-a_{\lambda,[1^2]}-a_{\lambda,[2^1]})(P_{8}- \lambda_1 P_{6}- (q^2+q+1-\lambda_1) P_{5}). \end{multline*}

\subsection{Weight $7$} \label{sec-w7}
If $P$ lies in $A_{\lambda,\emptyset}$ or $A_{\lambda,[1^1]}$ then the dimension of $L_P$ will not be the expected one. We will need another definition for the singular cases.

\begin{dfn} \label{dfn-Bb} Let $B_{\lambda,\mu}$ be the set of $\lambda$-tuples $P$ for which there are two distinct lines, both defined over $k$, whose intersection point $p$ lies in $P$, and such that one of the lines contains a $\mu$-tuple of points of $P\setminus p$ and the other a $\lambda-[1^1]-\mu$-tuple of points of $P\setminus p$. Define also
\begin{multline*} b_{\lambda,\mu}:= \binom{q^2+q+1}{2} \frac{\zee{\lambda}}{\zee{\mu} \zee{\lambda-[1^1]-\mu}}\abs{X(\mu)} \abs{X(\lambda-[1^1]-\mu)} \\ \text{where} \; X:=\Pbk^1 \setminus \Pbk^0. \end{multline*} 
\end{dfn}

\begin{lem} \label{lem-Bb} If $\abs{\mu} > 1$ and $\abs{\lambda-[1^1]- \mu} > 1$ then $\abs{B_{\lambda,\mu}}=b_{\lambda,\mu}$.
\end{lem}
\begin{pf} Similar to Lemma \ref{lem-Aa}. \end{pf}

The new special cases are when $P$ is in either $B_{\lambda,[1^3]}$, $B_{\lambda,[1^1,2^1]}$ or $B_{\lambda,[3^1]}$ and $S$ lies on any of the \emph{two} distinguished lines of $P$. Apart from this there is the case when $P$ lies in $A_{\lambda,[1^3]} \setminus B_{\lambda,[1^3]}$, $A_{\lambda,[1^1,2^1]} \setminus B_{\lambda,[1^1,2^1]}$ or $A_{\lambda,[3^1]} \setminus B_{\lambda,[3^1]}$ and then there is \emph{one} distinguished line of $P$ on which $S$ should lie. 

As an example, let us in detail show that if $P=(p,p_1,\ldots,p_6)$ and $S$ are not in any of the special cases described above, then $\dim L'_{P,S} = \dim L'_{P - p,S}-1$. For any distinct points $p$ and $q$, let $l(p,q)$ denote the line through $p$ and $q$. There are at most two points, say $p_1$ and $p_2$, that lie on the line through $p$ and $S$. There must also be two distinct points, say $p_3$ and $p_4$, such that the lines $l(p_1,p_3)$ and $l(p_2,p_4)$ do not contain $p$. Then the curve given by $l(p_1,p_3)$, $l(p_2,p_4)$, $l(S,p_5)$ and $l(S,p_6)$ is in $\dim L'_{P - p,S}$ but not in $L'_{P,S}$.

We thus find:
\begin{multline*}
\linsub{1}= a_{\lambda,\emptyset}(P_{9}- (q+1)P_{8}- q^2 P_{6}) +a_{\lambda,[1^1]}(P_{8}- (q+1)P_{7}-P_{6}- (q^2-1) P_{5})\\+a_{\lambda,[1^2]}(P_{7}- (q+1)P_{6}-2P_{5}- (q^2-2) P_{4})+a_{\lambda,[2^1]}(P_{7}- (q+1)P_{6}- q^2 P_{4}) \\
(b_{\lambda,[1^3]}+b_{\lambda,[1^1,2^1]}+b_{\lambda,[3^1]}) (P_{7}- (2q+1)P_{5}- (q^2-q) P_{4})\\
+ (a_{\lambda,[1^3]}-2b_{\lambda,[1^3]}) (P_{7}- (q+1)P_{5}- 3P_{5}-(q^2-3) P_{4})\\
+ (a_{\lambda,[1^1,2^1]}-2b_{\lambda,[1^1,2^1]}) (P_{7}- (q+1)P_{5}- P_{5}-(q^2-1) P_{4})\\
+ (a_{\lambda,[3^1]}-2b_{\lambda,[3^1]}) (P_{7}- (q+1)P_{5}- q^2 P_{4})\\
+(\Pbpl-a_{\lambda,\emptyset}-a_{\lambda,[1^1]}-a_{\lambda,[1^2]}-a_{\lambda,[2^1]}-a_{\lambda,[1^3]}-a_{\lambda,[1^1,2^1]}-a_{\lambda,[3^1]}\\+b_{\lambda,[1^3]}+b_{\lambda,[1^1,2^1]}+b_{\lambda,[3^1]}) \cdot (P_{7}- \lambda_1 P_{5}- (q^2+q+1-\lambda_1) P_{4}). \end{multline*}

\section{Cubic curves and elliptic curves} \label{sec-cubics}
In this section we will apply the method we have used for counting pointed non-singular quartics (described in Section \ref{sec-sieve}), to count pointed non-singular cubics. That is, we will choose $M=0$, apply the modified sieve method to the degree $3$ curves and then amend for the singular cubics. 

Let us reinterpret Definitions \ref{dfn-Lp} and \ref{dfn-linsub} where we let $L_{P}$ be the linear subspace of the $\Pb^9(k)$ of degree $3$ curves over $k$ that pass through $P$, and similarly Definition~\ref{dfn-t} where we also let the singular curves be of degree $3$. Define $\mathcal{C}(k)$ to be the set of non-singular cubic curves over $k$. What we will compute is 
\begin{equation} \label{eq-how-count-cub} 
\frac{1}{\PGLk{3}} \sum_{C \in \mathcal{C}(k)} \abs{C(\lambda)} = 
\frac{1}{\PGLk{3}} \bigl(\sum_{P \in \Pbk^2(\lambda)} \abs{L_P}-\sum_{\abs{\mu} > 0} t_{\lambda,\mu} \bigr), \end{equation}
for all partitions $\lambda$ with $\abs{\lambda} \leq 8$. 

In the first two subsections we will compute $t_{\lambda,\mu}$ for all partitions $\mu$ with $\abs{\mu} \geq 1$. Since the singular degree $3$ curves all have irreducible components whose normalizations have genus $0$, we will be able to do this for partitions $\lambda$ of any weight. In the next subsection $\abs{L_P}$ will be computed for all $[1^n]$-tuples $P$, where $n \leq 8$. In the final subsection this information is shown to be sufficient to compute equation \eqref{eq-how-count-cub} for any $\lambda$ with $\abs{\lambda} \leq 8$. We will then translate the counts of pointed non-singular cubics into $\s_n$-equivariant counts of points of the moduli space $\Mm{1}{n}$ for $n \leq 10$.

\subsection{The reducible degree $3$ curves} \label{sec-red-cubics}
All the reducible degree $3$ curves have components of degree $1$ or $2$ and they can be dealt with in the same way as the degree $4$ curves with such components. 

\begin{table}[htbp] \caption{Singular degree $3$ curves} \label{tab-cub}
\centerline{
\vbox{
\offinterlineskip
\hrule
\halign{&\vrule#& \quod #\hfil \quod \strut &\vrule#& \quod \hfil#\hfil \strut \quod \cr
height2pt&\omit&&\omit& \cr
& \bf{Description} && \bf{\#} & \cr
height2pt&\omit&&\omit&\cr
\noalign{\hrule}
height2pt&\omit&&\omit&\cr
&\emph{\small Triple line} && $\scriptstyle \frac{1}{q^3(q+1)(q-1)^2}$ &\cr
height2pt&\omit&&\omit&\cr
&\emph{\small Line + double line} &&$\scriptstyle \frac{1}{q^2(q-1)^2}$ &\cr
height2pt&\omit&&\omit&\cr
&\emph{\small $\lambda$-set of lines, int. in $\mu$-set} &&  &\cr
height2pt&\omit&&\omit&\cr
&$\scriptstyle \star \,\,\lambda,\mu = ([1^3],[1^3]),([1^1,2^1],[1^1,2^1]), ([3^1],[3^1])$ &&$\scriptstyle \frac{1}{6(q-1)^2},\frac{1}{2(q+1)(q-1)},\frac{1}{3(q^2+q+1)}$&\cr
height2pt&\omit&&\omit&\cr
&$\scriptstyle \star \,\,\lambda,\mu = ([1^3],[1^1]),([1^1,2^1],[1^1]), ([3^1],[1^1])$ &&$\scriptstyle \frac{1}{6q^2(q-1)},\frac{1}{2q^2(q-1)}, \frac{1}{3q^2(q-1)}$&\cr
height2pt&\omit&&\omit&\cr
&\emph{\small Conic + line, int. in $\lambda$-set} && &\cr
height2pt&\omit&&\omit&\cr
&$\scriptstyle \star \,\,\lambda= [1^1], [1^2], [2^1]$ &&$\scriptstyle \frac{1}{q(q-1)}, \frac{1}{2(q-1)}, \frac{1}{2(q+1)}$&\cr
height2pt&\omit&&\omit&\cr
} \hrule}}
\end{table}

\subsection{The singular cubics} \label{sec-cubsing}
The inverse image, under a normalization morphism, of the singularity of a singular cubic is a $\lambda$-set $P$ on the normalization $\Pbk^1$, where $\abs{\lambda}=1$ or $2$. Fix such a $\lambda$-set $P=\{p_1,p_2\}$ and consider the line $l$ of divisors of the form $p_1+p_2+r$ for any $r$ on $\Pbk^1$. Inside the $3$-dimensional space of all divisors of degree $3$ we let $\mathfrak{d}$ be any plane over $k$ that contains $l$, except for the two (possibly equal) planes (not necessarily defined over $k$) that are given by divisors of the form $p_1+r_1+r_2$ for any $r_1$, $r_2$ on $\Pbk^1$ and $p_2+s_1+s_2$ for any $s_1$, $s_2$ on $\Pbk^1$. This gives $q+1-\lambda_1$ possibilities to choose the linear system $\mathfrak{d}$. The two exceptional planes are the only ones that contain a base point. 

\begin{lem} \label{clm-cubsing} The image of any morphism induced by $\mathfrak{d}$ is a singular cubic. 
\end{lem}

\begin{lem} \label{lem-cubsing} The number of singular cubics, of a type determined by the partition $\lambda$, is equal to 
$$\frac{1}{\zee{\lambda}} \cdot \abs{\Pbk^1(\lambda)} \cdot  (q+1-\lambda_1) \cdot \frac{1}{\PGLk{2}}.$$
\end{lem}

\subsection{Linear subspaces of degree $3$ curves} \label{sec-linsub-cubics}
We will restrict ourselves to the case $\lambda=[1^{\abs{\lambda}}]$. In the next subsection we will see that this is the only case we need. Arguing as in Section \ref{sec-linsub} we get the following results.

\begin{itemize}
\item[$\star$] If $n \leq 4$ then $\linsub{0}= \abs{\Pbk^2([1^n])} P_{9-n}.$ 
\item[$\star$] If $n=5$ then $\linsub{0}= a_{[1^5],\emptyset}P_{5}+(\abs{\Pbk^2([1^5])}-a_{[1^5],\emptyset}) P_{4}. $
\item[$\star$] If $n=6$ then $\linsub{0}= a_{[1^6],\emptyset}P_{5}+a_{[1^6],[1^1]}P_{4}+(\abs{\Pbk^2([1^6])}-a_{[1^6],\emptyset}-a_{[1^6],[1^1]}) P_{3}. $
\item[$\star$] If $n=7$ then $\linsub{0}= a_{[1^7],\emptyset}P_{5}+a_{[1^7],[1^1]}P_{4}+a_{[1^7],[1^2]}P_{3} 
+(\abs{\Pbk^2([1^7])}-a_{[1^7],\emptyset}-a_{[1^7],[1^1]}-a_{[1^7],[1^2]}) P_{2}.$
\end{itemize}

When the weight is $8$ some new phenomena occur. If we have a pair of lines with four non-singular points over $k$ lying on each, then there is a $\Pb^2$ of degree $3$ curves passing through these points. The number of choices of such points is equal to $c:= 8! \binom{q^2+q+1}{2} \binom{q}{4}^2.$ By B\'ezout's Theorem, all degree~$3$ curves that contain seven points lying on a conic must contain that conic. Hence, if a set $P$ of eight points lie on a conic then $\dim L_P=2$. The number of choices of eight such points over $k$ is by Example \ref{exa-conic}, equal to $d:= 8! \cdot q^2 (q^3-1) \binom{q+1}{8}.$ The dimension of $L_P$ is the expected one in all other cases and we conclude that if $n=8$ then
\begin{multline*}
\linsub{0}= a_{[1^8],\emptyset}P_{5}+a_{[1^8],[1^1]}P_{4}+a_{[1^8],[1^2]}P_{3} 
+ a_{[1^8],[1^3]} P_{2}+ c P_{2} + d P_{2} \\
+(\abs{\Pbk^2([1^8])}-a_{[1^8],\emptyset}-a_{[1^8],[1^1]}-a_{[1^8],[1^2]}-a_{[1^8],[1^3]}-c-d) P_{1}. \end{multline*}

\subsection{From cubics to elliptic curves} \label{sec-elliptics}
For a non-singular irreducible projective curve $C$ over $k$, the Lefschetz trace formula tells us that $a_m(C)$ (defined in Section \ref{sec-firstpart}) is equal to $\mathrm{Tr}(F^m,H_c^1(C_{\bar k}, \Q_l))$ where $H_c^1(C_{\bar k}, \Q_l)$ denotes the compactly supported $\ell$-adic \'etale cohomology. There are $g$ pairs of eigenvalues $\alpha_i$, $\bar \alpha_i$ of the Frobenius $F$ acting on $H_c^1(C_{\bar k}, \Q_l)$ and thus $a_m(C)=\sum_{i=1}^g \alpha_i^m+\bar \alpha_i^m.$ These eigenvalues are also known to fulfill $\alpha_i \bar \alpha_i = q$.

If $g=1$, we can conclude that $a_m(C)$ is expressible in terms of $a_1(C)$ and $q$. In other words, if we can compute the number of points over $k$ of a non-singular genus~$1$ curve over $k$, we can compute the number of points over any finite extension of $k$. From equation \eqref{eq-m1n-count} it then follows that to make an $\s_n$-equivariant count of $\Mm{1}{n}$ it suffices to compute $\abs{\Mmt{1}{n}^{F \cdot id_n}}=\abs{\Mmt{1}{n}^{F}}$. 

\begin{lem} \label{lem-cub-ell} For any $n \geq 1$ we have 
$$\abs{\Mmt{1}{n+1}^{F}}+n \cdot \abs{\Mmt{1}{n}^{F}} = \frac{1}{\PGLk{3}} \sum_{C \in \mathcal{C}(k)} \abs{C(id_n)}.$$
\end{lem}
\begin{pf} Take a non-singular cubic over $k$ together with points $q_1,\ldots,q_n$ over $k$. The tangent line through $q_1$ either intersects in no other point or it intersects in a point $q$ over $k$, which is either distinct from or equal to one of the points $q_1,\ldots,q_n$. 

On the other hand, we can embed an element $(C,p_1,\ldots,p_{n+1})$ of $\Mmt{1}{n+1}^{F}$ with the complete linear system $|2p_1+p_{n+1}|$, and an element $(C,p_1,\ldots,p_{n})$ of $\Mmt{1}{n}^{F}$ in $n$ different ways using the complete linear systems $|2p_1+p_{i}|$ for $1 \leq i \leq n$. 
The image curves will be precisely the ones described above.
\end{pf}

In the previous subsections we have counted non-singular cubics with up to eight points over $k$ and in Example 7.1 of \cite{Jonas1} the well-known fact that $\abs{\Mmt{1}{1}^{F}}=q$ is shown. We can therefore use Lemma \ref{lem-cub-ell} to compute $\abs{\Mmt{1}{n}^{F}}$ for $n \leq 9$.

In \cite{Jonas1} it is shown that if $i$ is odd then
$$\sum_{[(C,p)] \in \Mm{1}{1}(k)/\cong_{k}} \frac{a_1(C)^i}{\abs{\Aut_k \bigl((C,p)\bigr)}} = 0.$$ 
Equation \eqref{eq-m1n-count} for $\sigma=id$ then shows that if $n$ is even we need only compute $\abs{\Mmt{1}{j}^{F}}$ for all $j \leq n-1$ to be able to conclude $\abs{\Mmt{1}{n}^{F}}$. In turn this shows that we can compute $\abs{\Mmt{1}{n}^{F}}$ for $n \leq 10$ and then by the opening arguments of this subsection also $\abs{\Mmt{1}{n}^{F \cdot \sigma}}$ for any $n \leq 10$ and permutation $\sigma \in \s_n$.

\bibliographystyle{plain}
\bibliography{cite}

\end{document}